
%
%
\input amstex
\input amsppt.sty

\magnification=\magstephalf
\CenteredTagsOnSplits
\NoBlackBoxes
\def\today{\ifcase\month\or
 January\or February\or March\or April\or May\or June\or
 July\or August\or September\or October\or November\or December\fi
 \space\number\day, \number\year}
\define\({\left(}
\define\){\right)}

\define\Aut{\operatorname{Aut}}
\define\CC{{\Bbb C}}

\define\Diff{\operatorname{Diff}}

\define\Hom{\operatorname{Hom}}
\define\Map{\operatorname{Map}}

\define\RR{{\Bbb R}}

\define\Tr{\operatorname{Tr}}
\define\ZZ{{\Bbb Z}}
\define\[{\left[}
\define\]{\right]}

\define\chiup{\raise.5ex\hbox{$\chi$}}
\define\cir{S^1}

\define\endexer{\bigskip\tenpoint}
\define\exertag #1#2{\removelastskip\bigskip\medskip\eightpoint\noindent%
\hbox{\rm\ignorespaces#2\unskip} #1.\ }

\define\inv{^{-1}}
\define\mstrut{^{\vphantom{1*\prime y}}}
\define\protag#1 #2{#2\ #1}

\define\res#1{\negmedspace\bigm|_{#1}}
\define\temsquare{\raise3.5pt\hbox{\boxed{ }}}
\define\theexertag{\theprotag}
\define\theprotag#1 #2{#2~#1}

\define\zmod#1{\ZZ/#1\ZZ}
\define\zn{\zmod{n}}

\def\bigstrut{\hbox{\vrule height18pt depth 9pt width0pt}}
\def\dashfill{\cleaders\hbox to .8em{\hss--\hss}\hfill}
\def\entry#1:#2:#3 {\bigstrut{#1}&{#2}&{#3}\cr}
\def\titlestrut{\hbox{\vrule height12pt depth 9pt width0pt}}
\define\C{\Cal{C}}
\define\E{\Cal{E}}

\define\G{\Cal{G}}
\define\LL{\Cal{L}}
\define\Pbar{\overline{P}}
\define\RZ{\RR/\ZZ}
\define\Rep{\operatorname{Rep}}
\define\Vect{\operatorname{Vect}}
\define\Xhat{\hat{X}}
\define\arr#1#2{#1@>{#2}>>}
\define\bX{\partial X}
\define\bY{\partial Y}
\define\catl{\Cal{L}}
\define\cut{^{\text{cut}}}
\define\eac#1#2{e^{\tpi S_{#1}(#2)}}
\define\fld#1{\Cal C_{#1}}
\define\fldb#1{\overline{\fld {#1}}}
\define\fldbp#1{\overline{\fld{#1}'}}
\define\form#1#2{\Omega ^#2_{#1}}
\define\hl#1#2{L_{#1}(#2)}
\define\hol{\operatorname{hol}}
\define\id{\operatorname{id}}
\define\lgh{\widehat{LG}}
\define\meas#1#2{d\mu _{#1}(#2)}
\define\mytimes{\odot}
\define\thg{\omega }
\define\tpi{2\pi i}
\define\zo{[0,1]}

\newread\epsffilein    
\newif\ifepsffileok    
\newif\ifepsfbbfound   
\newif\ifepsfverbose   
\newdimen\epsfxsize    
\newdimen\epsfysize    
\newdimen\epsftsize    
\newdimen\epsfrsize    
\newdimen\epsftmp      
\newdimen\pspoints     
\pspoints=1bp          
\epsfxsize=0pt         
\epsfysize=0pt         
\def\epsfbox#1{\global\def\epsfllx{72}\global\def\epsflly{72}%
   \global\def\epsfurx{540}\global\def\epsfury{720}%
   \ifx#1[\let\next=\epsfgetlitbb\else\let\next=\epsfnormal\fi\next{#1}}%
\def\epsfgetlitbb#1#2 #3 #4 #5]#6{\epsfgrab #2 #3 #4 #5 .\\%
   \epsfsetgraph{#6}}%
\def\epsfnormal#1{\epsfgetbb{#1}\epsfsetgraph{#1}}%
\def\epsfgetbb#1{%
%
%
\openin\epsffilein=#1
\ifeof\epsffilein\errmessage{I couldn't open #1, will ignore it}\else
%
%
   {\epsffileoktrue \chardef\other=12
    \def\do##1{\catcode`##1=\other}\dospecials \catcode`\ =10
    \loop
       \read\epsffilein to \epsffileline
       \ifeof\epsffilein\epsffileokfalse\else
%
%
          \expandafter\epsfaux\epsffileline:. \\%
       \fi
   \ifepsffileok\repeat
   \ifepsfbbfound\else
    \ifepsfverbose\message{No bounding box comment in #1; using defaults}\fi\fi
   }\closein\epsffilein\fi}%
%
%
\def\epsfsetgraph#1{%
   \epsfrsize=\epsfury\pspoints
   \advance\epsfrsize by-\epsflly\pspoints
   \epsftsize=\epsfurx\pspoints
   \advance\epsftsize by-\epsfllx\pspoints
%
%
   \epsfxsize\epsfsize\epsftsize\epsfrsize
   \ifnum\epsfxsize=0 \ifnum\epsfysize=0
      \epsfxsize=\epsftsize \epsfysize=\epsfrsize
%
%
     \else\epsftmp=\epsftsize \divide\epsftmp\epsfrsize
       \epsfxsize=\epsfysize \multiply\epsfxsize\epsftmp
       \multiply\epsftmp\epsfrsize \advance\epsftsize-\epsftmp
       \epsftmp=\epsfysize
       \loop \advance\epsftsize\epsftsize \divide\epsftmp 2
       \ifnum\epsftmp>0
          \ifnum\epsftsize<\epsfrsize\else
             \advance\epsftsize-\epsfrsize \advance\epsfxsize\epsftmp \fi
       \repeat
     \fi
   \else\epsftmp=\epsfrsize \divide\epsftmp\epsftsize
     \epsfysize=\epsfxsize \multiply\epsfysize\epsftmp
     \multiply\epsftmp\epsftsize \advance\epsfrsize-\epsftmp
     \epsftmp=\epsfxsize
     \loop \advance\epsfrsize\epsfrsize \divide\epsftmp 2
     \ifnum\epsftmp>0
        \ifnum\epsfrsize<\epsftsize\else
           \advance\epsfrsize-\epsftsize \advance\epsfysize\epsftmp \fi
     \repeat
   \fi
%
%
   \ifepsfverbose\message{#1: width=\the\epsfxsize, height=\the\epsfysize}\fi
   \epsftmp=10\epsfxsize \divide\epsftmp\pspoints
   \vbox to\epsfysize{\vfil\hbox to\epsfxsize{%
      \includegraphics{#1}%
      \hfil}}%
\epsfxsize=0pt\epsfysize=0pt}%

%
%
{\catcode`\%=12 \global\let\epsfpercent=
%
%
\long\def\epsfaux#1#2:#3\\{\ifx#1\epsfpercent
   \def\testit{#2}\ifx\testit\epsfbblit
      \epsfgrab #3 . . . \\%
      \epsffileokfalse
      \global\epsfbbfoundtrue
   \fi\else\ifx#1\par\else\epsffileokfalse\fi\fi}%
%
%
\def\epsfgrab #1 #2 #3 #4 #5\\{%
   \global\def\epsfllx{#1}\ifx\epsfllx\empty
      \epsfgrab #2 #3 #4 #5 .\\\else
   \global\def\epsflly{#2}%
   \global\def\epsfurx{#3}\global\def\epsfury{#4}\fi}%
%
%
\def\epsfsize#1#2{\epsfxsize}
%
%
\let\epsffile=\epsfbox

\refstyle{A}
\widestnumber\key{SSSSSSS}   
\document

	\topmatter
 \title\nofrills  Quantum Groups from Path Integrals \endtitle
 \author Daniel S. Freed  \endauthor
 \thanks The author is supported by NSF grant DMS-9307446, a Presidential
Young Investigators award DMS-9057144, and by the O'Donnell
Foundation.\endthanks
 \affil Department of Mathematics \\ University of Texas at Austin\endaffil
 \address Department of Mathematics, University of Texas, Austin, TX
78712\endaddress
 \email dafr\@math.utexas.edu \endemail
 \date January 25, 1995\enddate
 \subjclass 81R50, 81-01, 58Z05\endsubjclass
	\endtopmatter

\document

The goal of these lectures is to explain how quantum groups arise in
3~dimensional topological quantum field theories (TQFTs).  Of course,
``explain how'' is not the job of science, and perhaps you will find other
explanations more satisfying.  Let me explain!

What is a 3~dimensional TQFT?  At the very least it gives a topological
invariant of 3~dimensional manifolds.  That is, to each 3-manifold~$X$ it
assigns a complex number~$Z_X$ and if the invariants for~$X$ and~$X'$~are
different ($Z_X\not= Z_{X'}$), then $X$~and $X'$~are not diffeomorphic.  A
3~dimensional TQFT also gives invariants of knots and links.  For example, if
$K$~is a knot in ordinary 3-space, then we get a set of numerical
invariants~$\{I_K(\alpha )\}$ indexed by some finite set.  The {\it Jones
polynomial\/} of a knot, and similar polynomial invariants of knots, fit into
this picture.  As a historical note, Vaughan Jones~\cite{J} introduced his
polynomial invariant in the mid 80s before the advent to topological quantum
field theories.  Those were introduced by Edward Witten in 1987 (following a
suggestion of Michael Atiyah), first in 4~dimensions to give a quantum field
theoretic interpretation of Donaldson's invariants of 4-manifolds.  A few
years later~\cite{W} he introduced a 3~dimensional TQFT which reproduces the
Jones polynomial and which is our concern here.  The classical action of this
field theory is the {\it Chern-Simons invariant\/}, which was introduced into
geometry in the early 1970s.  As a mathematician I must immediately point out
that Witten's methods, involving the path integral, are far from an
established part of rigorous 1990s mathematics.

Shortly after Witten's paper on quantum Chern-Simons invariants, Reshetikhin
and Turaev~\cite{RT1} showed how to start with extremely complicated
algebraic data---called a {\it quantum group\/}---and again produce the Jones
polynomial and its generalizations.  (Their work is completely rigorous.)
Subsequently, they showed~\cite{RT2} how to use the same data to construct
invariants of 3-manifolds.  The construction of a complete TQFT from this
algebraic data, which involves more than invariants of 3-manifolds and knots,
has been folklore ever since and by now is completely written down~\cite{T}.
I remark here that instead of starting with a quantum group, one can start
with certain ``categorical'' data instead.

Now I can explain my point of view in these lectures.  The algebraic data of
either a quantum group or its categorical equivalent is extremely
complicated!  One could hardly guess in advance that such data can produce
invariants of knots and 3-manifolds.  Nor can one easily construct algebraic
data satisfying the necessary hypotheses.  By contrast the classical
Chern-Simons action is beautiful and simple!  It is relatively easy to write
down.  One sees from the beginning that Lie groups enter the picture in a
fundamental way.  And if you are willing to accept the path integral (you
shouldn't!), then you have a nice geometric construction of the Jones
polynomial and related 3-manifold invariants.  This leads us to pose the
following.

 \medskip
 \noindent{\bf Problem:} Start with the Chern-Simons action and construct the
quantum group which gives the same 3-manifold and knot invariants.
 \medskip

My goal is to explain how to do this in a simple case.  As I said, the
Chern-Simons theory starts with a compact Lie group~$G$ (and a piece of
topological data which I'll explain later).  The Jones polynomial concerns
the case $G=SU(n)$ for variable~$n$.  In the simple case we treat $G$~is a
{\it finite\/} group.  This was first considered by Dijkgraaf and
Witten~\cite{DW}.  The major simplification here is that the path integral is
a finite sum, rather than an integral over an infinite dimensional space, so
is rigorously defined.  So we immediately get a 3-manifold invariant, though
it is rather simple and relatively uninteresting.  The knot invariants are
possibly more interesting; I don't believe that they have been investigated
fully.  In any case our interest is in the quantum group and our strategy is
this: We exploit the fact that the path integral is well-defined to introduce
generalizations of the path integral.  Thus one ingredient in a 3~dimensional
TQFT is a ``quantum Hilbert space''~$E(Y)$ for every surface~$Y$.  In usual
quantum field theories it is constructed by canonical quantization.  In our
simple model we show how to get it by an exotic path integral.  Something is
immediately very strange---the result of an integration is a Hilbert space!
Even more strange is the path integral we introduce for a 1-manifold, i.e.,
for a circle~$S$.  There is where we will see the quantum group emerging.  In
fact, the quantum groups we compute this way were written down in a paper of
Dijkgraaf, Pasquier, and Roche~\cite{DPR}.  They did not related it to the
Chern-Simons invariant.  It was a conjecture of Altschuler and
Coste~\cite{AC} that these quantum groups construct (via the
Reshetikhin-Turaev prescription) the invariants of the finite group
Chern-Simons theory.  Our methods prove this conjecture.

This, then, is our strategy.  In a $d$~dimensional field theory, where
usually the classical action is only defined for fields in $d$~dimensions, we
will generalize the classical action to fields on manifolds of dimension less
than~$d$.  We then introduce a corresponding generalization of the path
integral for these exotic classical actions.  Of course, one is immediately
led to ask whether our constructions can be generalized, at least
heuristically, in Chern-Simons theory with continuous gauge group, or
possibly in other quantum field theories.  At this writing I do not know the
answer to this question.

The mathematics here is complicated and abstract, but not difficult.  I have
already tried my hand at several accounts.  The original paper, with all of
the computations, is~\cite{F1}.  Previous joint work with Frank
Quinn~\cite{FQ} discusses the basic theory in more detail.  The
brief~\cite{F2} gives a heuristic account of our extension of TQFT (without
mentioning path integrals or the classical theory) as well as a discussion of
central extensions (which arise in Chern-Simons theory with continuous gauge
groups) and invariants of ``framed tangles.''  The conference
proceedings~\cite{F3} contains a heuristic explanation of our generalized
path integrals as well as a brief idea of how some of this structure appears
in characteristic numbers in topology.  The summer school notes~\cite{F4} I
wrote a few years ago give a more leisurely introduction to the basics of
finite group Chern-Simons theory.  For these summer school notes I will
unabashedly cut and paste text and pictures from my previous writings!  Much
of the work I will leave for you in the form of guided exercises.  Many of
the exercises are not directly related to the topic of these lectures, but
perhaps you will find them useful anyhow.

\head
\S{1}: Classical Field Theory
\endhead
\comment
lasteqno 1@ 44
\endcomment

We begin with a discussion of the basic ingredients in a classical field
theory---spacetimes, fields, action.  This is mostly to fix the ideas and
notation since these same ingredients in the finite group Chern-Simons theory
may otherwise appear exotic.  We then discuss the Wess-Zumino-Witten action
in 2~dimensions in some detail.  Here we introduce an extension of the idea
of the classical action of a field theory.  Although this particular example
is not the subject of our lectures, it is a familiar example and hopefully
the geometry we discuss is easily accessible.

 \subhead Classical Actions
 \endsubhead

A field theory has a particular dimension attached to it, which we
call~`$d$'.  The standard examples of field theories have~$d=4$ and take
place on Minkowski space.  We allow more general examples.  This means first
of all that $d$~is not necessarily~4.\footnote{In the Chern-Simons example of
most interest to us~$d=3$.} Also, we allow {\it spacetimes\/} which are
curved manifolds.  Generally we take them to be compact, with or without
boundary.  Notice the terminology we will often use: A manifold is called
{\it closed\/} if it is compact and has no boundary.  Now Minkowski space
carries a metric whose isometry group is the Poincar\'e group, and one
usually requires that the the theory have the Poincar\'e group as a symmetry.
We will generalize this considerably.  Namely, a given field theory takes
place on manifolds with some specified extra structure.  This may be
topological (e.g. an orientation or spin structure) or may be geometric (a
conformal structure or a metric).  Thus usual {\it relativistic\/} field
theory should be thought of as taking place on Lorentz manifolds.  We allow
field theories on Riemannian manifolds and also on manifolds with less rigid
structure.  A field theory based on Riemannian manifolds is called a {\it
Euclidean field theory\/}, one based on manifolds with a conformal structure
is a {\it conformal field theory\/}, and one based on manifolds with only
some topological structure is a {\it topological quantum field theory\/}.  We
can be even more adventurous, of course.  We might allow singular manifolds,
for example.  For $d=1$ there are interesting field theories defined on
graphs.  Or we might allow supermanifolds.  Or we might take more abstract
sorts of spaces.

So the first ingredients of a classical field theory are a dimension and a
class of spacetimes,\footnote{We use the word `spacetimes' even though in our
context the words `space' and `time' may not have much significance.} i.e.,
manifolds of that dimension.  The next ingredient is a space of fields~$\fld
X$ for each spacetime.  This is usually a set of local functions on~$X$ if we
interpret `function' liberally enough.  The main object of study is the {\it
classical action\/} (or simply, action) which is usually a real-valued
function function on the space of fields:
  $$ S_X \: \fld X\longrightarrow \RR. \tag{1.1} $$
Typically, if a field is denoted~$\phi $, there is a lagrangian
density~$L_X(\phi)$ whose value at~$x\in X$ only depends a finite number of
derivatives of the field at~$x$ and on a finite number of derivatives of the
geometric data (such as a metric) on~$X$ at~$x$.  Then the action is the
integral of the Lagrangian density:
  $$ S_X(\phi ) = \int_{X}L_X(\phi). \tag{1.2} $$
The Lagrangian density is a density(!), that is, something which can be
integrated over~$X$.  If $X$~is oriented then it can be taken to be a
differential $d$-form.

Before reviewing the main properties of an action, let's note how some
familiar examples fit into this scheme.

        \exertag{1.3} {Exercise}
 Consider first classical mechanics as a field theory with $d=1$.  Look first
at a particle moving in~$\RR^3$.  Suppose that $X=[0,T]$.  Then $\fld X$~is
the space of paths in~$\RR^3$.  If there is no potential, then what is the
action?  What is the Lagrangian density?  What if there is now a potential
function $V\:\RR^3\to \RR$?  What if we replace~$\RR^3$ by an arbitrary
Riemannian manifold~$M$?
        \endexer

        \exertag{1.4} {Exercise}
 A more invariant formulation of classical mechanics is as a field theory on
1~dimensional Riemannian manifolds~$X$.  Can you formulate the Lagrangian
density and action in this case?  (The fields should be paths into a fixed
Riemannian manifold~$M$.)  Note that if $X$~is diffeomorphic to an interval,
then it is {\it isometric\/} to~$[0,T]$ for some~$T$.  In other words, the
only invariant of a Riemannian interval is its length.  This is a fancy way
of stating parametrization by arclength.
        \endexer

        \exertag{1.5} {Exercise}
 Generalize the previous example to the {\it $\sigma $-model\/} in
$d$~dimensions.  This is a $d$~dimensional field theory formulated on
Riemannian $d$-manifolds.  There is a fixed auxiliary Riemannian
manifold~$M$, and the fields on a $d$~dimensional Riemannian manifold~$X$ are
smooth maps $\phi \:X\to M$.  The Lagrangian density is
  $$ L_X(\phi ) = | d\phi |^2 d\mu _X,  $$
where $d\mu_X$ is the Riemannian volume density.  Write~$L_X$ in local
coordinates, or in a form recognizable to you.  Prove that this action is
{\it conformally invariant\/} in $d=2$~dimensions.
        \endexer

        \exertag{1.6} {Exercise}
 Here is a simple topological example in $d$~dimensions.  Fix a manifold~$M$
(no metric!) and a $d$-form $\omega \in \form Xd$.  The spacetimes are now
simply oriented $d$-manifolds~$X$ and the space of fields~$\fld X$ is the
space of smooth maps $\phi \:X\to M$ (as it is for the $\sigma $-model).  The
lagrangian density, or better lagrangian form, is $L_X(\phi )=\phi ^*\omega
$.  Write this in local coordinates.  Are there any simplifications if
$\omega $~ is closed ($d\omega =0$)?
         \endexer

        \exertag{1.7} {Exercise}
 Try to formulate free field theories in this formalism.  Consider for
example a free scalar field or a free spinor field.  What is the precise
class of spacetimes considered?  What is the space of fields?  Do you know
some interaction terms to add to these Lagrangians?
        \endexer

        \exertag{1.8} {Exercise}
 Consider now a gauge theory.  We will be doing this in more detail later,
but it is a good idea for you now to think of how this fits in with our
formalism.  To be concrete, consider Yang-Mills in 4~dimensions.  Now for a
Riemannian 4-manifold~$X$ the space of fields is a space of connections.  Can
you see how to fit gauge symmetry into the picture?  What is the action?
Does Yang-Mills make sense in other dimensions?  Show that only in $d=4$~is
the Yang-Mills lagrangian conformally invariant.  Show that in~$d=2$ it is
invariant under area-preserving diffeomorphisms.  (This is crucial in David
Gross' lectures.)
        \endexer

There are two main properties of fields and classical actions I want to
mention here: {\it symmetry\/} and {\it locality\/}.  I may use the term
`functoriality' for the symmetry I am talking about here, or perhaps
`external symmetry'.  It is {\it not\/} like gauge symmetry (\theexertag{1.8}
{Exercise}) which is a symmetry of the fields, something we might call an
`internal symmetry.'  Rather, it refers to the symmetries of the spacetimes.
Such a symmetry is a diffeomorphism $f \:X'\to X$ which preserves all of the
structure.  So if we are dealing with a Euclidean field theory, the map~$f $
is required to preserve the metrics, i.e., $f $~is an isometry.  Of course,
we might have~$X'=X$ which is the most interesting case.  In any case we
require that a symmetry induce a map on fields
  $$ f ^*\:\fld X\longrightarrow \fld{X'} \tag{1.9} $$
and that the action be preserved:
  $$ S_{X'}(f ^*\phi ) = S_X(\phi ),\qquad \phi \in \fld X. \tag{1.10}
     $$
If there is a Lagrangian density, then we usually have the stronger condition
that the density is preserved, at least up to an exact term.

\midinsert
\bigskip
\centerline{
 \epsfxsize=350pt
\epsffile{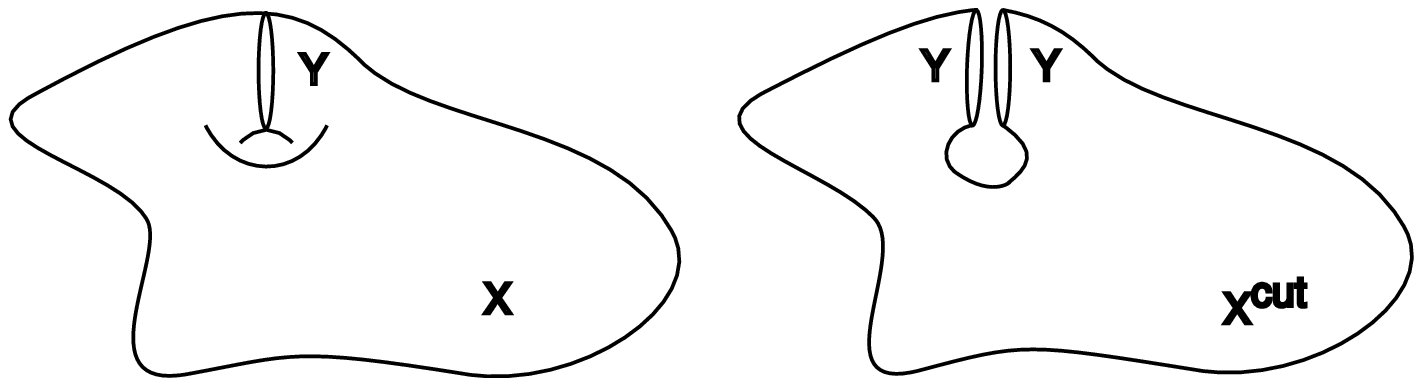}}
\nobreak
\medskip
\centerline{Figure~1: Cutting a spacetime~$X$ along~$Y$ to obtain~$X\cut$}
\bigskip
\endinsert

{\it Locality\/} is the assertion that the fields are local objects and that
the action can be computed locally.  There are many ways to formulate this,
and for our purposes we focus on the following situation.  Suppose $X$~is a
$d$~dimensional spacetime and $Y\hookrightarrow X$~is a closed codimension
one submanifold of~$X$.  If we cut~$X$ along~$Y$ then we obtain a new
spacetime which I will call~$X\cut$.  Notice that we do not require that our
spacetimes be connected, nor that they have connected boundaries.  The usual
picture is Figure~2, in which $X$~is connected and $X\cut$~has two
components.  But this is not necessary.  Nor is it necessary that $Y$~be
connected. In any case note that $X\cut$~has two new pieces in the boundary,
each of which is diffeomorphic to~$Y$.  (One of them appears with the
opposite orientation.)  The situation is illustrated in Figure~1.  Notice
that there is a gluing map $g\:X\cut\to X$ which identifies these two
boundary components.  Suppose now we have a field~$\phi $ on~$X$.  Then there
is a pullback field~$\phi \cut$ on~$X\cut$ and we require that
  $$ S_{X\cut}(\phi \cut) = S_X(\phi ). \tag{1.11} $$
This is trivial if the action is given by an integral as in~\thetag{1.2}.

\midinsert
\bigskip
\centerline{
 \epsfxsize=350pt
\epsffile{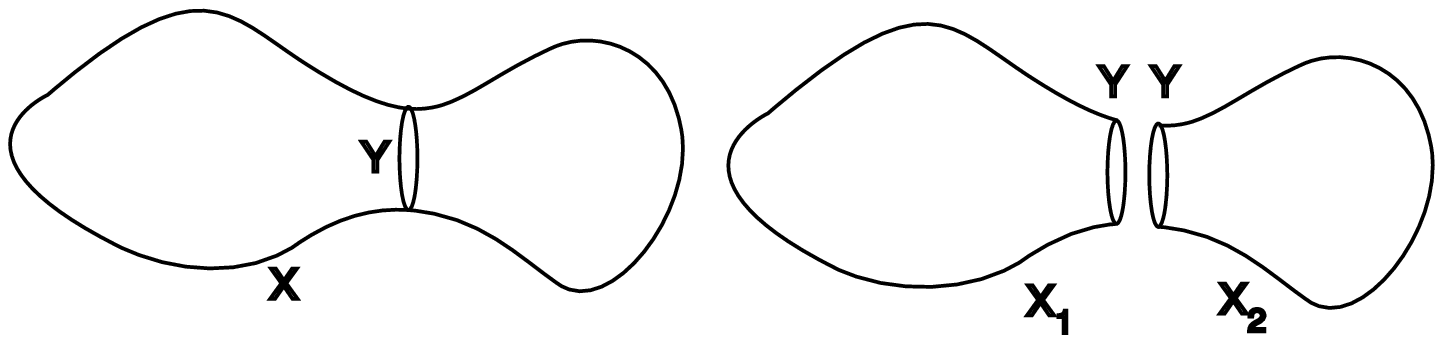}}
\nobreak
\medskip
\centerline{Figure~2: Cutting~$X$ into two pieces}
\bigskip
\endinsert

        \exertag{1.12} {Exercise}
 Consider a field theory formulated on Minkowski space.  Show that
\thetag{1.10}~is the assertion that the action is invariant under the
Poincar\'e group.  How does the Poincar\'e group act on scalar fields?  On
spinor fields?
        \endexer

        \exertag{1.13} {Exercise}
 Verify \thetag{1.10} and~\thetag{1.11} for the examples you considered
earlier.
        \endexer

        \exertag{1.14} {Exercise}
 If the spacetimes in a field theory have an orientation, then usually the
action satisfies:
  $$ S_{-X}(\phi )= -S_{X}(\phi ), \tag{1.15} $$
where $-X$~denotes the manifold~$X$ with the opposite orientation.  Check
this out in examples.  On field theories formulated on Minkowski space, what
does this say about Lorentz transformations not connected to the identity?
        \endexer

        \exertag{1.16} {Exercise}
 Can you think of more general formulations of locality?  What example(s)
motivates your formulation?
        \endexer

 \subhead The Wess-Zumino-Witten Action
 \endsubhead

As motivation for our extension of the notion of the classical action we
consider an action which is somewhat more interesting geometrically than
those indicated in the exercises above.  Often this {\it
Wess-Zumino-Witten\/} (WZW) action is just one term in the lagrangian.
In~$d=2$ it is added to the $\sigma $-model lagrangian (\theexertag{1.5}
{Exercise}) to define one of the fundamental conformally invariant field
theories in 2~dimensions.  From a mathematician's point of view the
corresponding quantum field theory teaches us much about the representation
theory of loop groups.  In any case, our interest for the moment is in the
geometry of the classical action.

Let $G=SU(2)$ be the Lie group of $2\times 2$~unitary matrices of determinant
one.  Recall that as a manifold $G$~is diffeomorphic to the 3-sphere.  Let
$g$~denote a general element of~$G$, and then
  $$ \theta = g\inv dg  $$
is a matrix-valued differential 1-form on~$G$.  We are interested in the
(scalar) differential 3-form
  $$ \thg = c\,\Tr(\theta ^3) \tag{1.17} $$
where $c$~is chosen so that
  $$ \int_{SU(2)}\thg = 1. \tag{1.18} $$

        \exertag{1.19} {Exercise}
 Parametrize~$SU(2)$ by matrices $\left(\smallmatrix \alpha &\beta
\\-\bar{\beta }&\bar{\alpha} \endsmallmatrix\right)$, where $\alpha ,\beta
\in \CC$ and $|\alpha|^2 + |\beta |^2 = 1$.  Write~$\theta $ and~$\thg$ in
terms of~$\alpha $ and~$\beta $.  Also, compute the constant~$c$.
        \endexer

        \exertag{1.20} {Exercise}
 More generally, let $G$~be any compact Lie group.  Define the analog of the
forms~$\theta $ and~$\thg$.  To do this you will need to pick an invariant
inner product on the Lie algebra of~$G$.  You can carry through the following
discussion for this more general case if $G$~is {\it simply connected\/}.  If
$G$~is not simply connected, then you will need much fancier (but
interesting) ideas.
        \endexer

The spacetimes in this $(1+1)$~dimensional field theory are oriented
2-manifolds (surfaces).  Let $X$~be a {\it closed\/} oriented surface.  Then
the fields we are interested in are simply maps to~$G$:
  $$ \fld X = \Map(X,G).  $$
Suppose that $\phi \:X\to G$ is such a map.  Since $G$~is diffeomorphic to a
3-sphere we can always extend to a map $\Phi \:W\to G$ where $W$~is an
oriented 3-manifold\footnote{Actually, it is better to say that $W$~is an
oriented 3-{\it chain\/} in the sense of homology theory.} whose boundary
is~$X$, and the restriction of~$\Phi $ to the boundary of~$W$ is~$\phi $.  We
write $\partial \Phi = \phi $.  Now define
  $$ S_X(\Phi )=\int_{W}\Phi ^*(\thg). \tag{1.21} $$
Provisionally, this is the action, but we want something which only depends
on~$\phi $, not on the extension~$\Phi $.  It is straightforward to see that
different choices of~$\Phi $ change~\thetag{1.21} by an integer (this
requires~\thetag{1.18}), so that
  $$ \eac X\phi := \eac X\Phi  \tag{1.22} $$
is independent of~$\Phi $.  This is the (exponentiated) WZW~action, and it
takes values in the complex numbers.  Note that since $\thg$~is real, it
always has unit norm.

        \exertag{1.23} {Exercise}
 Verify the assertion that \thetag{1.22}~is independent of the
extension~$\Phi $.
        \endexer

        \exertag{1.24} {Exercise}
 Prove that the WZW~action is functorial for orientation-preserving
diffeomorphisms $f \:X'\to X$, in the sense of~\thetag{1.9}
and~\thetag{1.10}.
        \endexer

        \exertag{1.25} {Exercise}
 Prove that the exponentiated WZW~action is a smooth function on~$\fld X$.
Compute its differential and find the critical points.
        \endexer

Notice that~$S_X(\phi )$ is {\it not\/} well-defined, only the
exponential~$\eac X\phi $ is.  This is fine from the path integral point of
view: Typically one writes the integrand of the path integral as $e^{i/\hbar\;
S_X(\phi )}$, and we can relax the normalization condition~\thetag{1.18} and
instead require that the ``coupling constant''~$c$ in~\thetag{1.17} only take
a discrete set of values which makes $e^{i/\hbar\; S_X(\phi )}$ well-defined.
This quantization of the coupling constant is well-known.  Notice also that
this action does {\it not\/} have a lagrangian density in the sense
of~\thetag{1.2}.

The interesting point is to define the WZW~action on an oriented surface~$X$
with nontrivial boundary~$\bX$.  Thus suppose $\phi \:X\to G$.  Since now
$\bX\not= \emptyset $ we can neither write~$X$ nor~$\phi $ as a boundary, and
so we cannot immediately find a 3-chain in~$G$ over which to integrate our
3-form~$\thg$.  Thus we proceed as follows.  The boundary~$\bX$ of~$X$ is a
disjoint union of circles.  Let $\Xhat=X\cup D$~denote the closed surface
obtained from~$X$ by gluing discs onto each of the boundary components.  The
union of those discs is denoted~$D$.\footnote{Notice that we allow~$\bX$ to
have multiple components.  Recall that we do not require that any of our
manifolds (for example~$X$) to be connected.} Extend~$\phi $ to a map~$\phi
\cup \Gamma \:\Xhat\to G$, where $\Gamma $~is the map on the union of the
discs.  Then the exponentiated action~$\eac{\Xhat}{\phi \cup\Gamma }$ is
well-defined.  However, we need to investigate its dependence on the
arbitrarily chosen~$\Gamma $.

        \exertag{1.26} {Exercise}
 Show that if $\Gamma '$~is another extension, then
  $$ \eac{\Xhat}{\phi \cup \Gamma '} = \eac{-D\cup D}{\Gamma \cup \Gamma
     '}\cdot \eac{\Xhat}{\phi \cup \Gamma } , \tag{1.27} $$
where $-D\cup D$ is the union of 2-spheres obtained by gluing $D$~with the
opposite orientation to another copy of~$D$.  (This is usually called the
{\it double\/} of~$D$.)  Now define
  $$ c_D(\Gamma ',\Gamma ) = \eac{-D\cup D}{\Gamma \cup\Gamma '}.
     \tag{1.28} $$
 Conclude that for maps $\Gamma ,\Gamma ',\Gamma ''\:D\to G$ we have
  $$ c_D(\Gamma '',\Gamma ) = c_D(\Gamma '',\Gamma ')\cdot c_D(\Gamma
     ',\Gamma ).   $$
        \endexer

Now we have a sort of exponentiated action for a field~$\phi $ on the
surface~$X$ with boundary.  It is not simply a number, however.  It is a
(complex-valued) function depending on how we extend~$\phi $ to~$\Xhat$, in
other words, depending on~$\Gamma $.  Furthermore, the dependence on~$\Gamma
$ is fairly simple~\thetag{1.27}, and most importantly the
factor~\thetag{1.28} which arises can be computed purely in terms of~$\Gamma
,\Gamma '$, i.e., it does not depend on~$\phi $.  (This is really a
manifestation of locality.)  One might be content to stop here and say that
this function of~$\Gamma $ {\it is\/} the exponentiated action.  However,
there is a nicer geometric way to proceed, and this is our reason for
considering this example.

Here, then, is our typical mathematician's ploy.  We have an example of a
function which obeys a certain equation~\thetag{1.27}.  So let's consider the
set of {\it all\/} functions which satisfy the same equation.  To do that we
first need to identify the set of allowable~$\Gamma $.  Let $\gamma =\partial
\phi $ be the restriction of~$\phi $ to the boundary~$\bX=\partial D$.  Then
the set of~$\Gamma $ is the set of fields on~$D$ whose restriction
to~$\partial D$ is~$\gamma $:
  $$ \fld D(\gamma ) = \{\Gamma \:D\to G : \partial \Gamma =\gamma \}.
      $$
Now we define the set of functions which satisfy~\thetag{1.27}, giving it the
suggestive name~$\hl{\bX}{\gamma} $:
  $$ \hl{\bX}{\gamma} = \{ \ell \:\fld D(\gamma )\to\CC : \ell (\Gamma ')=
     c_D(\Gamma ',\Gamma )\cdot \ell (\Gamma )\}. \tag{1.29} $$
Notice that this set depends only on~$\gamma $, as is indicated by the
notation.  Now what does $\hl{\bX}{\gamma}$~look like?  I claim that it is a
one dimensional complex vector space, also known as a {\it complex line\/}.
Furthermore, I claim that it has a natural inner product, so is actually a
{\it hermitian line\/}.

        \exertag{1.30} {Exercise}
 Prove these last two assertions by constructing the vector space structure
and inner product and verifying that $\hl{\bX}{\gamma} $~is one dimensional.
        \endexer

Finally, notice that the definition of~$\hl{\bX}{\gamma}$ does not use the
fact that $\bX$ is the boundary of a surface.  In other words, for any
field\footnote{You should take note here that although we originally
discussed fields in a $d$~dimensional field theory as defined on
$d$~dimensional manifolds (spacetimes), we are now extending that notion to
consider fields in a $d$~dimensional field theory defined on
$d-1$~dimensional manifolds.  Here it is clear that $\fld Y=\{\gamma \:Y\to
G\}$.} $\gamma \in \fld Y$ on a closed oriented 1-manifold\footnote{This is
simply a finite union of circles.}~$Y$ we use~\thetag{1.29} to define a
hermitian line~$\hl Y\gamma $.  (Then $D$~is the manifold obtained by
attaching a disk to each component of~$Y$.)

So to summarize we have used the 3-form~$\thg$ to define two mappings:
  $$ \align
      \phi \in \fld X &\longmapsto \eac X\phi \in \hl{\bX}{\partial \phi },
     \tag{1.31}\\
      \gamma \in \fld Y&\longmapsto \hl Y\gamma . \tag{1.32}\endalign $$
Here $X$ is a compact oriented 2-manifold and $Y$~is a {\it closed\/}
oriented 1-manifold.  In other words, we allow~$X$ to have a boundary, but
not~$Y$.  When $X$~is closed the line~$L_{\emptyset }$ is simply the
``trivial'' hermitian line of complex numbers~$\CC$.  The exponentiated
action is then defined by~\thetag{1.22}.  In case $X$~has boundary
\thetag{1.31}~is a convenient way to look at the construction preceding
\theexertag{1.26} {Exercise}.

Let me say it again: The exponentiated action on a manifold with boundary is
not a complex number, but now takes values in a hermitian line.  This is
already a modification of the scheme we outlined in the previous section, but
this action still has the essential locality property---or gluing
law---written in~\thetag{1.11}, only here it is expressed in a different
form---equation~\thetag{1.34} of the following exercise.

        \exertag{1.33} {Exercise}
 Let $X$~be a compact oriented 2-manifold, $Y\hookrightarrow X$ a closed
codimension one submanifold, and $X\cut$~the manifold obtained by cutting~$X$
along~$Y$ (Figure~1).  Let $\phi \in \fld {X}$ and $\phi \cut\in \fld{X\cut}$
the corresponding field on~$X\cut$, as in~\thetag{1.11}.  Prove that
  $$ \Tr_Y \eac{X\cut}{\phi \cut} = \eac X\phi , \tag{1.34} $$
where $\Tr_Y$~is performed using the hermitian metric in the line~$\hl Y{\phi
\res Y}$.  You will need to use the properties~\thetag{1.36}
and~\thetag{1.37} below to make sense of~\thetag{1.34}.
        \endexer

Now we make an even larger extension/modification of the scheme in the
previous section, and this is our key point in this section.  {\it We
consider \thetag{1.32}~as the definition of an (exponentiated) WZW~action for
fields on a 1-manifold.\/} This is not at all the usual picture.  First of
all, we have an action in a $d$~dimensional theory defined for fields in
$d-1$~dimensions.  Secondly, the value of the action is not a {\it number\/},
but rather a {\it set\/} (more precisely, a hermitian line).  I hope that
through the course of these lectures you will be convinced that this is a
useful extension of the notion of a classical action.  For now, here are some
properties which are analogous to properties of the usual classical action.
In particular, \theexertag{1.38} {Exercise} deals with symmetry.  Locality
will have to wait for the next lecture.  (Think about what this would mean.)

        \exertag{1.35} {Exercise}
 Construct isomorphisms
  $$ \hl{Y_1\sqcup Y_2}{\gamma _1\sqcup \gamma _2}\longrightarrow
     \hl{Y_1}{\gamma _1} \otimes \hl{Y_2}{\gamma _2} \tag{1.36} $$
and
  $$ \hl{-Y}\gamma \longrightarrow \overline{\hl Y\gamma }. \tag{1.37} $$
Here `$\sqcup $'~denotes the disjoint union (that is, the union of two sets
with empty intersection) and $\overline{L}$~is the complex conjugate vector
space to~$L$ (that is, the vector space with the same underlying addition and
the complex conjugate scalar multiplication).  Notice that these properties
are analogous to properties of the usual classical action.  (The
isomorphism~\thetag{1.37} is analogous to the equation~\thetag{1.15}.  We did
not write the analog of~\thetag{1.36} for the usual action, but you should
easily see what it is.  In fact, you can view it as a special case of the
gluing law where we glue along an empty manifold!) Show that these
isomorphisms are actually {\it isometries\/}.
        \endexer

        \exertag{1.38} {Exercise}
 Suppose $f \:Y'\to Y$ is an orientation-preserving diffeomorphism of
1-manifolds.  Then for any~$\gamma \in \fld Y$, construct an isometry
  $$ f ^*\: \hl Y\gamma \longrightarrow \hl{Y'}{f ^*\gamma }.  \tag{1.39} $$
Notice that~\thetag{1.39} is the analog of~\thetag{1.10}.
        \endexer

The next exercise is worked out in~\cite{F5}.  It is a geometric way to pass
from the lagrangian picture (classical action) to the hamiltonian picture.

        \exertag{1.40} {Exercise}
 Show that $\hl Y\gamma  $ depends smoothly on~$\gamma $ in the sense that
these lines fit together into a smooth hermitian line bundle $L_Y\to \fld Y$.
Define a parallel transport on paths in~$\fld Y$ using the WZW~action.  Show
that this is actually the parallel transport of a connection on~$L_Y$.
Compute the curvature of this connection.
        \endexer

        \exertag{1.41} {Exercise}
 You can carry out the constructions in this section in much more generality.
One generalization is to let $G$~be any compact Lie group, though you will
have to work harder if $G$~is not simply connected.  Another generalization
is to work in $d$~dimensions and replace~$\thg\in \Omega ^3(G)$ by a
$(d+1)$-form on a manifold~$M$.  The normalization condition~\thetag{1.18}
should be replaced by the condition that the integral of this form over all
$(d+1)$-cycles is an integer.  You will also want to make some topological
assumptions on~$M$ generalizing simple connectivity in the $d=2$~case.  So a
simpler case is a $d=1$~field theory based on an integral 2-form~$\omega $ on
some manifold~$M$.  In this case we can find a hermitian line bundle $L\to M$
with a unitary connection whose curvature is~$2\pi i\omega $.  Show that the
exponentiated classical action in this case is the parallel transport (or
holonomy) of this connection, and the action~\thetag{1.32} just reproduces
the line bundle~$L$.  (This last assertion is not as precise as it might
be---what is the precise statement?)
        \endexer

Suppose that $Y=\cir$ is the standard circle.  Then
$\fld{\cir}=\Map(S^1,G)=LG$ is the {\it loop group\/}; the multiplication of
loops is defined pointwise using the multiplication in~$G$.  According to
\theexertag{1.40} {Exercise} the lines~$L(\gamma )=\hl {\cir}\gamma $ fit
together to form a smooth hermitian line bundle $L\to LG$ over the loop
group.  Let $\lgh$~denote the set of elements of unit norm in~$L$; it is a
principal circle bundle over~$LG$.  In the next exercise you will show that
$\lgh$~is a central extension of~$LG$.  This construction of the central
extension is originally due to J. Mickelsson~\cite{M}.

        \exertag{1.42} {Exercise}
 Suppose $\gamma_1, \gamma _2\in LG$.  Construct an isometry
  $$ L(\gamma _1)\otimes L(\gamma _2)\longrightarrow L(\gamma _1\gamma _2).
     \tag{1.43} $$
This is not trivial---it uses a formula sometimes attributed in the physics
literature to Polyakov.  Restrict~\thetag{1.43} to the elements of unit norm
to define multiplication in~$\lgh$.  Verify that this multiplication is
associative and indeed defines a group.  Construct a homomorphism $\lgh\to
LG$ and show that its kernel is isomorphic to the circle group of unit
complex numbers.  Also, show that the kernel is {\it central\/} in~$\lgh$,
that is, elements in the kernel commute with every element in~$\lgh$.
        \endexer

        \exertag{1.44} {Exercise}
 \theexertag{1.38} {Exercise} asserts that the action of~$\Diff\cir$ on~$LG$
lifts to the line bundle~$L\to LG$.  How does the lifted action interact with
the isometry~\thetag{1.43}?
        \endexer

For details on the material in this section, see~\cite{F5}.

\head
\S{2}: Categories, Finite Groups, and Covering Spaces
\endhead
\comment
lasteqno 2@ 34
\endcomment

We discuss in general terms what ideas are necessary to extend the notion of
classical action further.  The indicated extension is crucial in
understanding the relationship of quantum groups to 3~dimensional TQFT.
Unfortunately it involves the concept of a {\it category\/}, which may be
off-putting at first.  Category theory has been called ``the theory of
abstract nonsense.''  Be that as it may, the notion is useful to us here.  We
then introduce what is surely the simplest field theory: gauge theory with
finite gauge group.  The interesting structure is in the space of fields; the
action we consider is trivial.  In the exercises we indicate a ``twisted''
version of the theory which has nontrivial action.  This theory exists in any
dimension, though our main interest later is in the $d=3$~case.  Then it is a
simple example of {\it Chern-Simons theory\/}, which can be defined for any
compact gauge group.  We remark that although these theories are analytically
simple, they still illustrate some basic properties of gauge theory,
especially the role of symmetries and reducible connections.

 \subhead Going Further
 \endsubhead

What I hope you learned from the WZW example is the following.  In a
$d$~dimensional field theory we can allow actions which are not of the
form~\thetag{1.1} but rather of the form
  $$ \eac X\cdot \:\fld X\longrightarrow \CC, \tag{2.1} $$
where $S_X(\cdot )$~may not be defined.  Furthermore, I challenge you to
think of any loss (from a physics point of view) in passing from~\thetag{1.1}
to the exponentiated form~\thetag{2.1} when $S_X(\cdot )$~{\it is\/} defined.
(I cannot think of any.)  Here $X$~is a {\it closed\/} $d$-manifold.  We
extend the idea of fields and action to closed $(d-1)$-manifolds~$Y$:
  $$ \hl Y\cdot \:\fld Y\longrightarrow \catl, \tag{2.2} $$
where $\catl$~is the {\it category\/} of all finite dimensional Hilbert
spaces.  We will have more to say about that shortly, so please don't panic
yet!  Finally, if $X$~is a $d$-manifold with boundary, the exponentiated
action~\thetag{2.1} has a generalization which we can explain using the
following diagram of line bundles:
  $$ \CD
      r^*L_{\bX} @>>> L_{\bX}\\
      @VVV @VVV\\
      \fld X @>r>> \fld{\bX}\endCD \tag{2.3} $$
Here $r$~is the restriction map which restricts a field to the boundary.  The
line bundle $L_{\bX}\to \fld{\bX}$ is the extended action~\thetag{2.2}.  Then
the action~$\eac X\cdot $ is a section of the line bundle $r^*L_{\bX}\to \fld
X$.  This extended classical action has several properties, scattered in
Lecture~1, which basically capture the idea that the action behaves like the
integral over the manifold of something which depends locally on the field.
The most characteristic of these properties is the {\it gluing
law\/}~\thetag{1.34}.

        \exertag{2.4} {Exercise}
 Perhaps an analogy with honest integration will help understand this idea of
an extended action.  The usual situation is that we have a compact oriented
$d$-manifold~$X$ and a $d$-form~$\alpha $ on~$X$.  Then the integral
$\int_{X}\alpha $ is defined and is a single number.  More generally,
consider a fiber bundle $X\to \Cal{C}$ whose typical fiber is a compact
oriented $d$-manifold, and let $\alpha \in \Omega ^d(X)$ be a d-form.  Then
{\it integration over the fibers\/} of $X\to Y$ produces a function
$\int_{X/Y}\alpha $ on~$\Cal{C}$.  This is the analog of the usual classical
action.  But more generally suppose that the fibers of $X\to\Cal{C}$ are
compact oriented manifolds of dimension~$d-i$.  Then integration along the
fibers gives an $i$-form on the base:
  $$ \int_{X/Y}\alpha \in \Omega ^i(\Cal{C}).  $$
So here we have a family of $(d-i)$-manifolds and we integrate something
$d$~dimensional over the fibers to get something $i$~dimensional on the base.
If $i=1$ then we get a 1-form, which you might think of as analogous to a
connection form on a hermitian line bundle.  This is analogous to our first
extension of the classical action, which in this case is not a 1-form but
rather a hermitian line bundle.

Incidentally, the exercise here is to fill in the details if you are not
already familiar with integration over the fibers.  You might also want to
consider the analog of Stokes' theorem in this context.
        \endexer

You should now raise several questions.  First, how does our extended
classical action work in familiar examples?  Secondly, are there nontrivial
examples other than the WZW action?  Finally, can we go further and consider
$(d-1)$-manifolds with boundary, closed $(d-2)$-manifolds, etc.?  Following
the old joke, we will answer these questions in the form of questions!
First, the familiar examples.

        \exertag{2.5} {Exercise}
 For the usual examples you considered in the previous section, show that the
extended action~\thetag{1.32} is trivial.  However, for usual second order
lagrangians the correct space of fields on the boundary should include a
derivative.  Consider classical mechanics, for example (\theexertag{1.4}
{Exercise}).  The space of fields attached to a point~$Y=\text{pt} $ should
be the tangent bundle~$TM$, not the manifold~$M$.  Then the extended action
gives a trivial line bundle over~$TM$.  However the connection constructed
following the idea of \theexertag{1.40} {Exercise} is nontrivial, and its
curvature is the standard symplectic form constructed from the Riemannian
metric.  Work out the details to the point that you recognize familiar
formulas from classical mechanics.
        \endexer

There is another nontrivial example: the Chern-Simons action.  I know this
best for~$d=3$, but in principle it can be worked out in other (odd)
dimensions.  (The $d=1$~case is \theexertag{1.41} {Exercise}.)  Our main
interest in these lectures is the Chern-Simons action for a finite gauge
group.  The case of a continuous gauge group has more geometric interest, and
you can find all of the details in~\cite{F5}.  (There is also an account
in~\cite{F4}.)  We consider the finite gauge group case in the next section
and defer to these references for the continuous group case.

Consider again a $d$~dimensional field theory.  We defined (in an example) an
action on $d$-manifolds with boundary which satisfies a gluing
law~\thetag{1.34}.  Further, we asserted that \thetag{1.32}~should be
considered as an extension of the classical action to closed
$(d-1)$-manifolds.  Now we want to go further---define an action on
$(d-1)$-manifolds with boundary and formulate a gluing law.  Let's just see
what kind of objects we should expect to run into.  Since the action on a
closed $(d-1)$-manifold gives a hermitian line, we expect that on a
$(d-1)$-manifold with boundary the action is some similar object.  (The
analogy in $d$~dimensions is that the action on a $d$-manifold with
boundary---an element in a complex line---is similar to the action on a
closed $d$-manifold---a complex number.)  At the very least we expect that it
is a set rather than some kind of number.  Then the analog of
equation~\thetag{1.34}, the gluing law, will be an ``equation'' between sets.
Now such equations are possible---you can say that two sets are equal---but
it is also possible to say that two sets are {\it isomorphic\/} without being
equal.  This is an extra layer of complexity which sets have that numbers
don't.  We have already run into this in \thetag{1.36}, \thetag{1.37},
\thetag{1.39}, and~\thetag{1.43}.  Another way to put it is that sets have an
``internal structure'' and it is possible to have {\it automorphisms\/} of
this structure.  Mathematicians (notably Saunders MacLane~\cite{Mc}) have
systematized these ideas in the notion of a {\it category\/}.  It fits into
the progression:
  $$ \text{number,\qquad set,\qquad category.}   $$

A category is a collection of {\it objects\/} and {\it morphisms\/} (maps and
arrows) between objects.  Two morphisms compose if the second begins where
the first ends.  The composition is assumed associative, and usually one
assumes that there are identity morphisms as well.  If we focus on the
objects, then the morphisms encode the internal structure that they possess.
On the other hand, it is useful to focus on the morphisms as well.  For
example, a category with one object is simply a set (of morphisms) with an
associative composition law and an identity, also known as a {\it
semigroup\/}.  If a category has more than one object, then we can think of
it as a ``semigroup with states.''  The objects represent the states, and
in each state there are certain morphisms which are possible.  Some change
the state and others (automorphisms) do not.

        \exertag{2.6} {Exercise}
 As an example consider the category~$\Cal{V}$ of finite dimensional complex
vector spaces and linear maps.  So an object in~$\Cal{V}$ is a vector space
and a morphism $L\:V_1\to V_2$ between $V_1,V_2\in \Cal{V}$ is a linear map.
Note that two vector spaces can be isomorphic but not equal.  Also, every
vector space of positive dimension has nontrivial automorphisms.  You are
familiar with the idea that we should consider isomorphic, but distinct,
vector spaces as being different.  Think of the tangent spaces to the
standard 2-sphere.  Tangent spaces at different points are isomorphic, but if
we could truly think of them as equal (in a ``continuous'' way) we would
quickly construct an everywhere nonzero vector field on the 2-sphere, which
is not possible.  Do you know other examples of categories?  Other situations
in which it is not possible to identify isomorphic objects which are not
equal?  What if two objects are isomorphic and have no nontrivial
automorphisms?  Can we safely identify them in that case?
        \endexer

        \exertag{2.7} {Exercise}
 Next we consider maps of categories.  (They are usually called {\it
functors\/}.)  Roughly, a functor maps the objects and morphisms of one
category into the objects and morphisms of another so that it preserves
compositions.  As an example, consider the functor $\Cal{V}\to\Cal{V}$ which
assigns to each vector space~$V\in \Cal{V}$ its double dual~$V^{**}$.  What
does this do to morphisms?
        \endexer

        \exertag{2.8} {Exercise}
 The extra layer of structure in a category allows us to define maps between
functors, called {\it natural transformations\/}.  Suppose that
$\Cal{C}_1,\Cal{C}_2$ are categories and
$\Cal{F}_1,\Cal{F}_2\:\Cal{C}_1\to\Cal{C}_2$ are functors.  Then a natural
transformation $\theta \: \Cal{F}_1\to\Cal{F}_2$ is for each object $C\in
\Cal{C}_1$ a morphism $\theta (C)\:\Cal{F}_1(C)\to\Cal{F}_2(C)$ such that it
is compatible with morphisms.  This means that for any morphism $C'@>f>>C$
in~$\Cal{C}_1$ the diagram
  $$ \CD
      \Cal{F}_1(C') @>\Cal{F}_1(f)>> \Cal{F}_1(C)\\
      @V\theta (C')VV @VV\theta (C)V\\
      \Cal{F}_2(C') @>\Cal{F}_2(f)>> \Cal{F}_2(C)\endCD  $$
commutes.  Construct a natural transformation from the identity functor to
the double dual functor of \theexertag{2.7} {Exercise}.  Show that this is
in fact a natural {\it isomorphism\/}.
        \endexer

Now we can make an educated guess about going further.  The action of a field
on a closed $(d-2)$-manifold should be a category, the action of a field on a
$(d-1)$-manifold with boundary should take values in the category associated
to the boundary field, and there should be a gluing law analogous
to~\thetag{1.34} which is a morphism in a category.  We can be even more
precise.  The type of category where the action takes its values should fit
into the progression:
  $$ \text{complex number of unit norm,\qquad hermitian line,\qquad ?.}
     \tag{2.9} $$
Before describing what the `?'~is, it is easier to consider the progression
of {\it trivial\/} values for the action:
  $$ \text{$1\in \CC$,\qquad $\CC$,\qquad $\catl$.}  \tag{2.10} $$
In these lectures we only consider a trivial action, so an understanding
of~\thetag{2.10} will suffice.  Recall that $\catl$~is the category of all
finite dimensional Hilbert spaces.  A morphism in~$\catl$ is an isometry.
This category is analogous to the complex numbers as the following exercise
shows.

        \exertag{2.11} {Exercise}
 We construct a structure on~$\catl$ analogous to the {\it ring\/} structure
on~$\CC$, that is, the addition and the multiplication.  There is a further
structure, which is complex conjugation.  Then we can build a norm from
multiplication and complex conjugation.  Construct, then, an addition functor
$\catl\times \catl\to \catl$, a multiplication functor $\catl\times \catl\to
\catl$, and a complex conjugation functor $\catl\to\catl$.  You should use
direct sum, tensor product, and the conjugate linear space.  Show that these
functors satisfy desired properties, but be careful that where there are
equalities in the complex numbers (associativity, commutativity, identity
element) there are natural isomorphisms in~$\catl$.  For example, the
statement that $\CC\in \catl$ acts as a multiplicative identity is the
assertion that the functor
  $$ \aligned
      \catl&\longrightarrow \catl\\
      L&\longmapsto \CC\otimes L\endaligned  $$
is naturally isomorphic to the identity functor.  So part of the structure
on~$\catl$ is the explicit specification of this natural isomorphism.  There
is a further layer of structure: {\it equations\/} among these natural
isomorphisms.  For example, associativity of addition is expressed by a
natural isomorphism
  $$ \theta _{v_1,V_2,V_3}\: (V_1\otimes V_2)\otimes V_3\longrightarrow
     V_1\otimes (V_2\otimes V_3).  $$
Now for four Hilbert spaces $V_1,V_2,V_3,V_4$ there is an equation among the
various~$\theta $'s.  What is this equation?
        \endexer

In the main flow of these lectures we will not need the general case of~`?'
in~\thetag{2.9}, only the trivial case in~\thetag{2.10}.  I don't mean to
make the~`?' so mysterious, and I'll give a nontrivial example in the
exercises.  Perhaps you have already realized that `?'~fits into the analogy
  $$ \multline\text{complex numbers} :\text{a hermitian line} \\=
     \text{category of finite dimensionsal Hilbert spaces} : \text{?}
     \endmultline\tag{2.12} $$
We will later introduce the idea of a ``2-Hilbert space'' and we will see
that a~`?' is a one dimensional 2-Hilbert space.  In any case we now
postulate that continuing the progression in~\thetag{2.1} and~\thetag{2.2} we
have for a {\it closed\/} $(d-2)$-manifold~$S$ an assignment
  $$ \Cal{L}_Y(\cdot) \:\fld S\longrightarrow \text{collection of one
     dimensional 2-Hilbert spaces}, \tag{2.13} $$
which satisfies properties analogous to those satisfied by~\thetag{2.1}
and~\thetag{2.2}.

        \exertag{2.14} {Exercise}
 Here is a category~$\Cal{K}$ which is a nontrivial example of a one
dimensional 2-Hilbert space.  Namely, let~$\Cal{K}$ be the category of all
unitary representations of the unitary group~$U(n)$ which are isomorphic to a
direct sum of several copies of the determinant representation on~$\CC$,
where the action of a unitary matrix is multiplication by the determinant.  A
morphism in~$\Cal{K}$ is an isomorphism of the representations (it commutes
with the action of~$U(n)$.)  Construct a functor $\catl\times
\Cal{K}\to\Cal{K}$ using the tensor product.  This is analogous to scalar
multiplication.  (Here $\catl$~are the scalars.)  Construct also the analog
of vector addition, namely a functor $\Cal{K}\times \Cal{K}\to\Cal{K}$, using
the direct sum.  Show that any irreducible (one dimensional) representation
$K_0\in \Cal{K}$ is a ``basis'' of~$\Cal{K}$.  In other words, use~$K_0$ and
scalar multiplication to construct an isomorphism $\catl\cong \Cal{K}$.
Consider how this fits with the analogy~\thetag{2.12}.  Can you provide a
rough definition for a 2-Hilbert space?
        \endexer

        \exertag{2.15} {Exercise}
 State explicitly the properties of~\thetag{2.13} which are ``analogous to
those satisfied by~\thetag{2.1} and~\thetag{2.2}.''  Note in particular the
gluing law for~\thetag{2.2}:
  $$ L_Y(\gamma )\cong \bigl(\hl{Y_1}{\gamma _1},\hl{Y_2}{\gamma _2}
     \bigr)_{\Cal{L}_S(\partial \gamma )}. \tag{2.16} $$
        \endexer

        \exertag{2.17} {Exercise}
 Here is an example of a higher dimensional 2-Hilbert space.  It will be
important in Lecture~4.  Let $G$~be a finite group and let $\Cal{E}$~denote
the category of all finite dimensional unitary representations of~$G$.
Construct an addition $\E\times \E\to\E$ and a scalar multiplication
$\catl\times \E\to\E$.  Then construct an inner product $\E\times
\E\to\catl$.  Can you find an ``orthonormal basis'' of~$\E$?  What is~$\dim
\E$?  This particular example of a 2-Hilbert space is also an algebra.  So
construct a multiplication $\E\times \E\to\E$.
        \endexer

        \exertag{2.18} {Exercise}
 Here is an example of a category which is ``small'' by comparison with the
previous examples.  It should dispel any illusion you may have that
categories are huge.  It is also an example which will be crucial for us
later.  Consider a finite group~$G$.  Construct a category whose objects are
elements of~$G$.  For each pair of elements~$x,g\in G$ we postulate a
morphism $g\:x\to gxg\inv $.  What is composition?  Show that if~$|G|=n$ this
constructs a category with $n$~objects and $n^2$~morphisms.  Draw a picture
for an abelian group.  For~$G=S_3$.
        \endexer

 \subhead Finite Gauge Theory
 \endsubhead

Fix a {\it finite\/} group~$G$. For example, $G$ could be the cyclic group
$\zn$ of $n$ elements, or the symmetric group $S_n$ of $n!$ elements.  There
is no restriction on $G$.  We need no other data to define our field theory,
though more data is needed for the twisted theory.

        \exertag{2.19} {Exercise}
 For the twisted theory we need to use the {\it classifying space\/}~$BG$.
One model is the following.  Construct an embedding $G\hookrightarrow O(N)$,
where $N=|G|$~is the order of~$G$.  Let $EG$~denote the space of all
$N$-tuples of orthonormal vectors in an infinite dimensional dimensional real
Hilbert space.  Show that $EG$~is contractible and that $O(N)$~acts freely
on~$EG$.  Define $BG$~to the the quotient $EG/G$, where $G$~acts via its
embedding into~$O(N)$.  The space $BG$~has interesting topology, and in the
twisted theory we need to fix a cohomology class\footnote{More precisely, we
need to fix a cocycle representing this cohomology class in some model of
cohomology.} $\lambda \in H^d(BG; \RR/\ZZ)$.  The untwisted theory in the
text corresponds to~$\lambda =0$.
        \endexer

The theory we consider is based on ``bare'' manifolds.  No orientation,
metric, etc. is needed.  Thus a spacetime is simply a compact $d$-manifold~$X$
with no additional structure.  (The twisted theory is based on oriented
manifolds, however.)  Symmetries of the spacetimes are simply diffeomorphisms.

Next we define the space of fields on a spacetime~$X$.  Recall that in our
extended notion of field theory we also consider fields on manifolds of
dimension less than~$d$, so we will define a space of fields for {\it any\/}
manifold~$M$.  Here is the definition:
  $$ \fld M = \left\{ \matrix P\\ \downarrow\\ M\endmatrix : \matrix P \text{
     is a {\it principal\/} ({\it Galois, regular\/})}\\ \text{ {\it covering
     space\/} with structure group } G.\endmatrix\right\} \tag{2.20} $$
In other terms, a field is a principal bundle $P\to M$ with structure group
$G$.  Thus $P$~is a manifold, the group $G$ acts freely on $P$, and the
quotient is $P/G = M$. We always take $G$ to act on the {\it right\/}.

As an example, consider $G= \zmod3$. If $M= \text{pt} $, then any bundle
looks like
  $$ \matrix \vdots &P \\ \downarrow&\\ \cdot& X\endmatrix $$
with $G$ cyclically permuting the 3 points. If $M=S^1$ then there are 3
possibilities, up to isomorphism, as illustrated in Figure~3.  The nontrivial
coverings are pictured as pieces of helices, but the endpoints are meant to
be identified.  Topologically, the total space~$P$ in these covers is a
circle.  The total space of the trivial cover is the disjoint union of
3~circles.

\midinsert
\bigskip
\centerline{
 \epsfysize= 2in
\epsffile{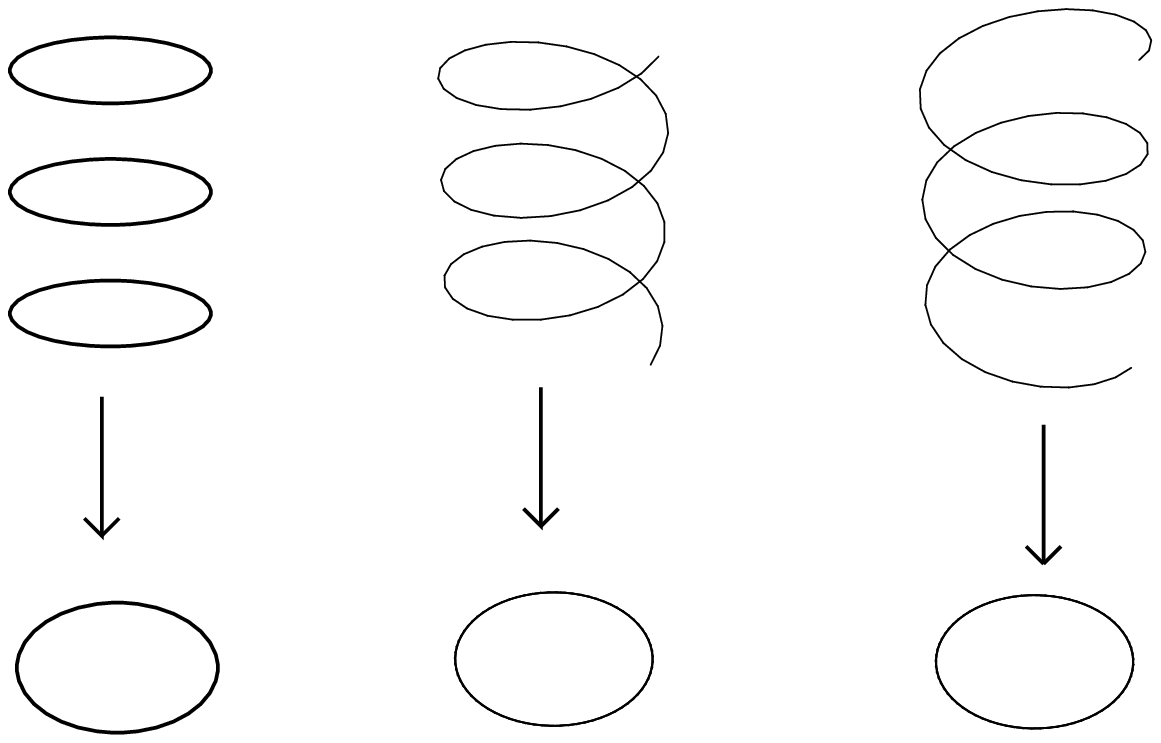}}
\nobreak
\medskip
\centerline{Figure~3: Principal $\zmod3$ bundles over~$\cir$}
\bigskip
\endinsert

You should have noticed that the space of fields~\thetag{2.20} is
qualitatively different than in the examples in Lecture~1.  In usual
examples the fields for a smooth manifold (usually infinite dimensional).
Here do they not even technically form a set!  (The collection of all of
anything is not a set---remember Russell!)  Rather, they form a category.
{\bf Warning:\/} Categories enter here in a different way than in our
extended notion of classical action in Lecture~1.  Remember that categories
are a fundamental mathematical structure, just as sets are, and we should not
be surprised to see them in a variety of contexts.  (We shall meet them again
in yet another context later.)  What the idea of a category captures here is
{\it gauge symmetry\/}.  These are symmetries of the fields (internal
symmetries), not symmetries of the spacetimes (external symmetries).  They
are the morphisms in the category of fields~$\fld M$.  Since they are all
invertible, we call them `isomorphisms'.

        \proclaim{\protag{2.21} {Definition}}
 An {\it isomorphism\/} $\varphi$ from $P'\in \fld M$ to $P\in \fld M$ is a
diffeomorphism $\varphi \: P'\to P$ which commutes with the $G$~action and
such that the induced map on the quotient $M$ is the identity.
        \endproclaim

\noindent
 This means that for $p'\in P'$ and $g\in G$ we have $\varphi (p'\cdot g) =
\varphi (p')\cdot g$, where the first~`$\cdot $' indicates the $G$~action
on~$P'$ and the second~`$\cdot $' the $G$~action on~$P$.  Any such map
induces a map $\overline{\varphi}\:M\to M$, and we restrict our isomorphism
to have $\overline{\varphi}=\id$.  More generally, we can consider diagrams
  $$ \matrix P'&\buildrel \varphi\over \to &P\\ \downarrow&&\downarrow\\
    M'&\buildrel {\bar\varphi}\over \to &M\endmatrix $$
where $\bar\varphi $~is possibly nontrivial.  They enter when we consider
symmetries of the manifolds as well as symmetries of the fields.

        \exertag{2.22} {Exercise}
 Let $f\:M'\to M$ be a diffeomorphism.  Construct $f^*\:\fld M\to\fld {M'}$
as required by~\thetag{1.9}.  In this case it is a functor.
        \endexer

The previous exercise shows how to pull back fields under maps of spacetimes.
For these to properly be considered fields we must also be able to cut and
paste, as in Figure~1.

        \exertag{2.23} {Exercise}
 Suppose $M$~is a manifold and $N\hookrightarrow M$~is a closed codimension
one submanifold of~$M$.  Let $M\cut$~be the manifold obtained by cutting~$M$
along~$N$.  Let $g\:M\cut\to M$ be the gluing map.  Construct a functor
  $$ g^*\:\fld M\longrightarrow \fld{M\cut}.  $$
        \endexer

\midinsert
\bigskip
\centerline{
 \epsfysize=2 in
\epsffile{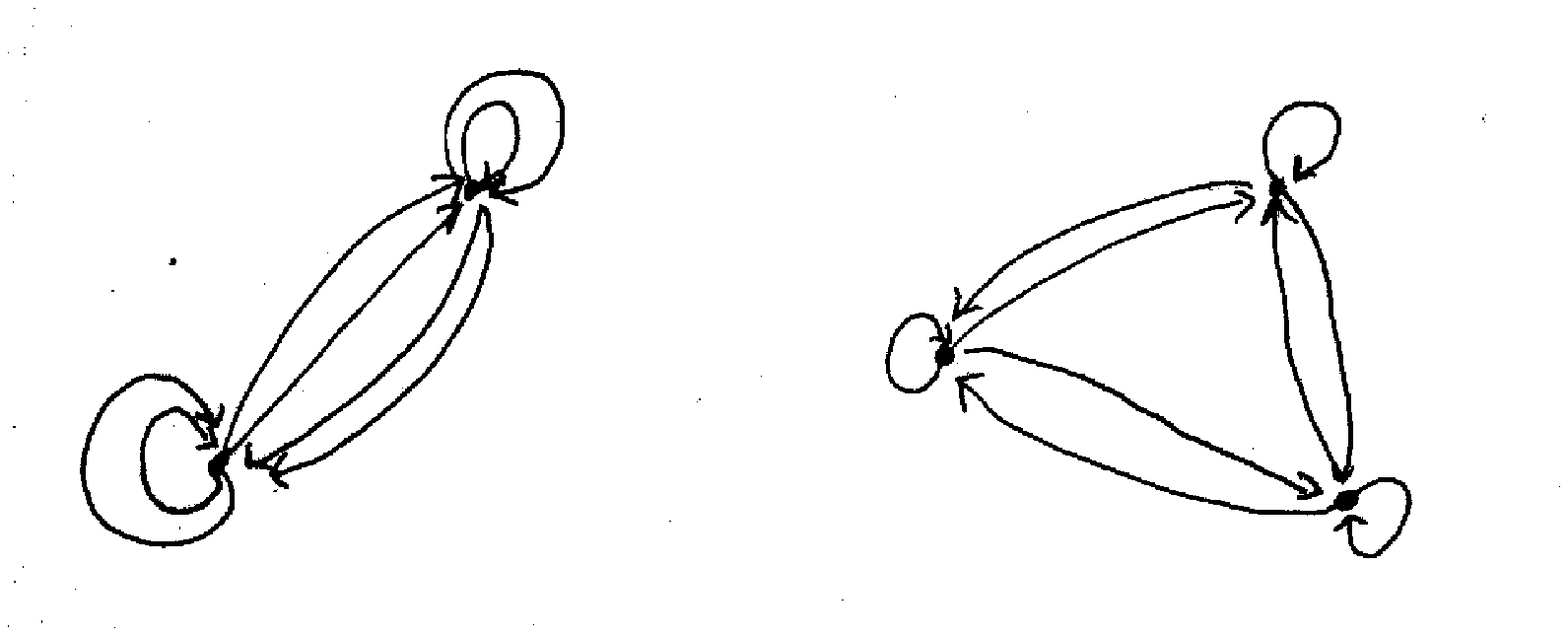}}
\nobreak
\medskip
\centerline{Figure~4:  The space of fields~$\fld M$}
\bigskip
\endinsert

So we get a picture of $\fld M$ which we schematically render in Figure~4.
Each point is a field and each arrow is a symmetry of fields, that is, an
isomorphism.  The arrows which start and end at an object~$P$ form a {\it
group\/} called $\Aut P$, the automorphism group of~$P$.  There is a finite
number of arrows from any point to any other.  $\Aut P$ is also called the
group of {\it gauge transformations\/} of~$P$.  We remark that this picture
applies to the space of (gauge) fields in any gauge theory, except that in
the general case there are continuous parameters---the space of connections
on a fixed bundle---as well as discrete ones---the choice of the bundle.
See~\cite{F5,\S1} for a discussion.

The fields $P$ and $P'$ are {\it equivalent\/} $(P\cong P')$ or {\it
isomorphic\/} if there is an arrow between them. Let
  $$ \fldb{M} = \hbox{ set of equivalence classes of fields on } M.
      $$
What makes gauge theories with finite gauge group tractable is that $\fldb
M$~is a finite set if $M$~is compact.

        \exertag{2.24} {Exercise}
 Determine $\fldb{M}$ for $M=pt$. For $M=S^1$.  For $M=[0,1]$.
        \endexer

        \exertag{2.25} {Exercise}
 Show that $\fldb M$ is a finite set for any compact manifold~$M$.
        \endexer

        \exertag{2.26} {Exercise}
 I claim that we can {\it not\/} make a field theory where the space of fields
associated to~$M$ is~$\fldb M$.  This is because we cannot paste equivalence
classes.  Consider, for example, cutting the circle~$\cir$ into an
interval~$\zo$.  Show that we can cut equivalence classes of bundles, but we
cannot paste them.
	\endexer

We can determine the space of equivalence classes of fields in terms of the
{\it fundamental group\/}.  Suppose $M$ is connected. Fix {\it basepoints\/}
$m\in M$ and $p\in P_m$, where $P_m$ is the fiber of $P$ at $m$.  Then a
field $ P\to M$ determines a map
  $$ \{\hbox{loops in~$M$ based at }m\}\longrightarrow  G$$
by taking the {\it holonomy\/} around the loop using the basepoint $p$.  (See
Figure~5.)  Any loop at~$m$ lifts uniquely to a {\it path\/} in~$P$ which
starts at~$p$ and ends at some~$p'$ in the fiber~$P_m$ of~$P$ over~$m$.  The
holonomy is the unique $h\in G$ satisfying $p'=p\cdot h$.  The holonomy only
depends on the homotopy class of the loop.

\midinsert
\bigskip
\centerline{
 \epsfysize=2 in
\epsffile{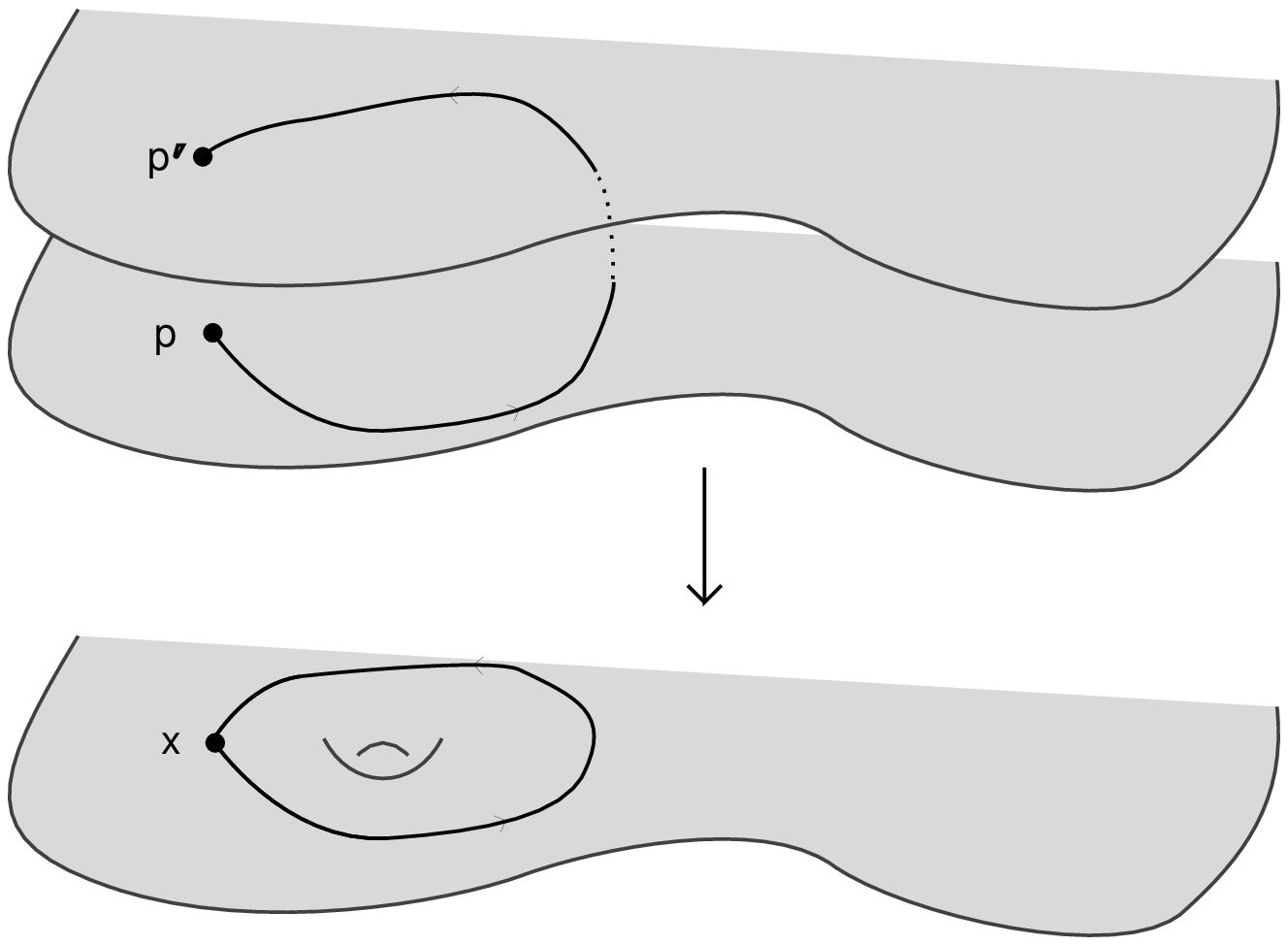}}
\nobreak
\medskip
\centerline{Figure~5: Definition of holonomy}
\bigskip
\endinsert

        \exertag{2.27} {Exercise}
 \item{a)} Check this last assertion.
 \item{b)} Show that the map $\pi_1 (M,m)\to G$ defined by holonomy is a
homomorphism of groups.
 \item{c)} If $\gamma$ is this homomorphism, and we change the basepoint $p$
and/or $m$, then the new homomorphism is $g\gamma g^{-1}$ for some $g\in G$.
 \item{d)} Suppose $M$~is connected and $h\:\pi _1(M,m)\to G$ is given.
Construct a bundle~$P\to M$ and a basepoint~$p\in P_m$ such that the holonomy
is~$h$.
	\endexer

\flushpar
 The assertions in \theexertag{2.27} {Exercise} hold more generally for any
{\it flat connection\/} with arbitrary gauge group.  In this case they show
that there is an isomorphism (of sets)
  $$ \fldb M \cong \Hom  (\pi_1(M),G)/G  \tag{2.28} $$
if $M$ is connected. (We omit the irrelevant basepoint in the notation for
the fundamental group.)

Finally we are ready to define the action in this theory.  We can summarize
in one word: The action is {\it trivial\/}!  Recall that we consider a
$d$~dimensional theory, so begin with the usual (numerical) action on
$d$-manifolds~$X$:
    $$ S_X (P) =0\ \hbox{ for all }\ P \in \fld X.  \tag{2.29} $$
So this is an example where the usual action \thetag{1.1}~is defined, not
just the exponentiated action~\thetag{2.1}.  The extended
actions~\thetag{2.2} and~\thetag{2.13} are also trivial.  So for a closed
$(d-1)$-manifold~$Y$ we have
  $$ L_Y(Q) =\CC\ \hbox{ for all }\ Q \in \fld Y, \tag{2.30} $$
and for a closed $(d-2)$-manifold~$S$ we have
  $$ \Cal{L}_S(R)=\catl\ \hbox{ for all }\ R \in \fld S. \tag{2.31} $$
Definitions~\thetag{2.29} and~\thetag{2.30} also make sense for manifolds
with boundary.

        \exertag{2.32} {Exercise}
 Verify the properties (discussed in Lecture~1) of the extended
action---gluing laws, functoriality, etc.  The functoriality statement is
more involved since one has to worry now about internal symmetries.  What is
the right formulation?
        \endexer

        \exertag{2.33} {Exercise}
 In the twisted theory the action is nontrivial~\cite{FQ}, ~\cite{F1}.  We
give a brief indication here.  Recall that the twisted theory is defined by
fixing a cohomology class~$\lambda \in H^d(BG;\RZ)$, or more precisely a
cocycle representing the cohomology class.  Prove first that for any
$G$~bundle $P\to M$ here exists a ``classifying map'' $F\:P\to EG$ which
commutes with the $G$~action, and that the map is unique up to homotopy
through $G$-maps.  So the quotient map $\overline{F}\:M\to EG$ is unique up
to homotopy.  If now $P\in \fld X$ for $X$~a closed {\it oriented\/}
$d$-manifold, we define
  $$ \eac XP = e^{\tpi \overline{F}^*(\lambda )[X]}, \tag{2.34} $$
for any classifying map~$F$.  Here $[X]\in H_d(X)$ is the fundamental class
of~$X$ given by the orientation.  Verify that \thetag{2.34}~is well-defined.

The twisted version of~\thetag{2.30} is more complicated.  Try to arrive at
it much as we arrived at~\thetag{1.29} in our study of the WZW~action.
Namely consider first the expression~\thetag{2.34} for a compact oriented
$d$-manifold~$X$ with boundary.  Now the orientation class is a relative
homology class---$[X]\in H_d(X,\bX)$---and the evaluation
$\overline{F}^*(\lambda )[X]$ does not make sense.  Instead you must choose
$d$-cycles on~$X$ which represent~$[X]$ and evaluate the chosen $d$-cocycle
(which
represents~$\lambda $) on such cycles.  This depends nontrivially on the
choice of the $d$-cycle, and the dependence is encoded in a complex line which
just depends on the restriction of~$P$ to~$\bX$.  Fill in the details and you
will have the twisted version of~\thetag{2.30}.

If you are truly ambitious you will carry this further and construct the
twisted version of~\thetag{2.31}.  This is another motivating example for our
extension of the classical action.

Try to formulate (if not prove) the properties of the extended action for
this twisted theory.  The symmetry properties are perhaps more apparent here
since the action is nontrivial.  The $d$~dimensional action
(twisted~\thetag{2.29}) is invariant under gauge symmetry, but what is the
analogous statement for the $d-1$~dimensional action (twisted~\thetag{2.30})?
What about for the $d-2$~dimensional action (twisted~\thetag{2.31})?
        \endexer

\head
\S{3}: Generalized Path Integrals
\endhead
\comment
lasteqno 3@ 34
\endcomment

Now we move from the classical theory to the quantum theory.  I dare say we
are following Feynman, though we are so far from the original context of his
path integrals that I suspect he would be either amused or more likely
appalled!  Our case is much simpler than usual examples in that the path
integral reduces to a finite sum.  Nonetheless, the formal picture---rigorous
here---is the same in usual examples.  We simplify our initial discussion by
ignoring the gauge symmetry.  The key new idea is an extension of the path
integral.  We interpret the usual quantum Hilbert space as the result of a
path integral over fields on a $(d-1)$-manifold (in a $d$~dimensional
theory).  Naturally it involves the extended classical action, as does our
extension to fields on $(d-2)$-manifolds.  We work out the usual path
integral and quantum Hilbert space in the finite group gauge theory.  In
Lecture~4 we will use the extended path integral to compute the quantum group
relevant to those theories.

 \subhead Path Integral Quantization
 \endsubhead

Imagine a $d$~dimensional field theory in which the space of fields~$\fld M$
attached to any manifold~$M$ (of dimension~$\le d$) is a finite set.  You may
as well assume that this is a topological field theory.  To simplify the
picture even further we assume that there is no (gauge) symmetry among the
fields.  We put back in the symmetry later.  However, here we allow a
nontrivial action.  The {\it partition function\/} of a closed
$d$-manifold~$X$ is defined to be
  $$ Z_X = \int_{\fld X} \eac XP\; \meas XP. \tag{3.1} $$
We write a typical field as~`$P$', keeping in mind our example.  The new
ingredient is a measure~$d\mu _X$ on the space of fields~$\fld X$.  With our
assumption that $\fld X$~is a finite set, the measure is simply a positive
number~$\meas XP$ for each~$P\in \fld X$.\footnote{Of course, in more typical
examples $\fld X$~is an infinite dimensional space and a measure~$d\mu _X$ is
extremely difficult to construct and has not been constructed in many
examples of interest.  Nonetheless, workers in {\it constructive quantum
field theory\/} have enjoyed many nontrivial successes in this pursuit.} Note
that the partition function~$Z_X$ is a complex number.

        \exertag{3.2} {Exercise}
 Prove that $Z_X$~is a topological invariant.  That is, if $f\:X'\to X$ is a
diffeomorphism, then $Z_{X'}=Z_X$.  Remember, we are assuming that this is a
topological field theory.
        \endexer

The preceding exercise is the invariance of the partition function under
symmetry.  It is the quantum analog of~\thetag{1.10}.  We now want to
investigate the quantum analog of the gluing law~\thetag{1.11}.  The
following discussion is the usual argument for locality of the path integral
via factorization into intermediate states.  In our language it goes as
follows.  Let $X$~be a $d$-manifold, which for simplicity you may assume to
be closed, and $Y\hookrightarrow X$ an embedded closed $(d-1)$-manifold.  We
assume that cutting~$X$ along~$Y$ splits~$X$ into two manifolds~$X_1$
and~$X_2$ (see Figure~2).  We identify $\partial X_1 = \partial X_2 = Y$ The
fields fit into the following diagram:
  $$ \CD
     \fld X @>c>> \fld{X_1}\times \fld{X_2}\\
     @Vr_1VV @VVr_2V\\
     \fld Y @>\Delta >> \fld Y\times \fld{Y}
     \endCD  $$
The vertical arrow~$r_1$ is restriction to~$Y$, the arrow~$r_2$ is
restriction to the boundaries of the~$ X_i$, the arrow~$\Delta $ is the
diagonal inclusion, and $c$~is the pullback under the gluing map.
(`$c$'~stands for `cutting'.)  Then we propose to do the integral over~$\fld
X$ in two stages using Fubini's theorem: First integrate over the fibers
of~$r_1$ and then over~$\fld Y$.  Now the gluing law for the classical
action~\thetag{1.34} says that if $P\in \fld X$ and $\langle P_1,P_2
\rangle\in \fld{X_1}\times \fld{X_2}$ the cut field, then
  $$ \eac XP = \(\eac{X_1}{P_1},\eac{X_2}{P_2}\)_{L_Y(Q)} \tag{3.3} $$
where $Q$~is the restriction of~$P$ to~$Y$.  The right hand side is the inner
product in the line~$\hl YQ$.  Carrying out the integration using the Fubini
theorem we obtain:
  $$ \split
      Z_X&=\int_{\fld X}\eac XP \;d\mu _X(P)\\
      &= \int_{\fld Y}\int_{r_1\inv (Q)} \eac XP \;d\mu _{r_1\inv (Q)}(P) \;
     d\mu _Y(Q)\\
      &= \int_{\fld Y} \int_{r_2\inv (Q,Q)} \eac XP \;d\mu _{r_2\inv
     (Q,Q)}(P\cut) \; d\mu _Y(Q)\\
      &= \int_{\fld Y}\(\int_{\fld {X_1}(Q)} \eac{X_1}{P_1}\;d\mu
     _{X_1}(P_1)\,,\,\int_{\fld {X_2}(Q)} \eac{X_2}{P_2}\;d\mu
     _{X_2}(P_2)\)_{\!\!\!\!\hl YQ}\!\!\!\!\!\!\!\!\!\!\!d\mu _Y(Q).\endsplit
     \tag{3.4} $$
Here we use the definition
  $$ \fld {X_i}(Q)=\{P\in \fld {X_i} : \partial P=Q\},\qquad Q\in \fld Y.
      $$
Also, we make certain implicit compatibility assumptions about the measure to
make this computation.

        \exertag{3.5} {Exercise}
 State explicitly these compatibility assumptions.
        \endexer

To rewrite this last expression in a nicer form, we make the following
definitions. Let~$X'$ be any $d$-manifold with boundary.
Then the path integral is a function of a field~$Q$ on the boundary:
  $$ Z_{X'}(Q) = \int_{\fld {X'}(Q)} \eac {X'}P \;\meas {X'}P,\qquad Q\in
     \fld{\partial X'}. \tag{3.6} $$
Note that the integral is over the space of fields with fixed boundary
value~$Q$.  The right hand side of~\thetag{3.6} takes values in the hermitian
line~$\hl {\partial X'}Q$, which is the extended action of the field on the
boundary (cf.~\thetag{2.3}).  So $Z_{X}$~is a section of the hermitian line
bundle $L_{\bX}\to \fld{\bX}$.  Again we use the mathematician's ploy of
introducing the space of {\it all\/} such sections.  Hence for any closed
surface~$Y$ set
  $$ E(Y) = \text{space of sections of the hermitian line bundle $L_{Y} \to
     \fld{Y}$.} \tag{3.7} $$
Note that $E(Y)$~is a finite dimensional complex vector space.  Then
\thetag{3.6} determines a {\it relative\/} invariant
  $$ Z_{X'} \in E(\partial X').  $$
We impose an $L^2$~inner product on~\thetag{3.7} using a measure~$d\mu _Y$
on~$\fld Y$ and the metric on the hermitian line bundle $L_{Y}\to \fld{Y}$.
With these definitions we rewrite~\thetag{3.4} as:
  $$ Z_{X} = \( Z_{X_1},Z_{X_2}\)_{E(Y)}. \tag{3.8} $$
This is the quantum gluing law---the quantum analog of~\thetag{3.3}---the
statement that the path integral is local.

Thus we have the standard formal ingredients of quantum field theory---the
partition function~$Z_X$ on closed spacetimes and the quantum Hilbert
space~$E(Y)$ for closed ``spaces'', i.e., $d$-manifolds.  We have indicated
that symmetry and locality hold for the partition function.  You should now
prove that other formal properties of the path integral and quantum Hilbert
space follow from the corresponding properties of the classical action.
Again you will have to make some assumptions about the measures~$d\mu _X$
and~$d\mu _Y$, which I leave you to formulate.

        \exertag{3.9} {Exercise}
 (Functoriality)\ Show that a diffeomorphism $f\:Y'\to Y$ induces an isometry
$f_*\:E(Y')\to E(Y)$.  In a field theory formulated on Minkowski space this
is a representation of the Euclidean group of a space slice on the quantum
Hilbert space; the representation of the whole Poincar\'e group involves the
path integral as well.  Show that a diffeomorphism $F\:X'\to X$ preserves the
partition function.  (Consider the case where $X$~has nontrivial boundary.)
        \endexer

        \exertag{3.10} {Exercise}
 (Multiplicativity)\ What can you say about~$E(Y_1\sqcup Y_2)$?  What about
$Z_{X_1\sqcup X_2}$?
        \endexer

        \exertag{3.11} {Exercise}
 Try this point of view out in some familiar field theories, even those
formulated in Minkowski space.  In that case, what does the gluing law say
about the usual propagation?  (Here you should start by considering ordinary
quantum mechanics, which is the case~$d=1$.)  A good nontrivial example is
the $d=2$ theory considered in David Gross' lectures: QCD~in 2~dimensions
with a fixed gauge group.  What is the Hilbert space~$E(S^1)$ in this theory?
Note that this is not a topological theory, but the partition function
depends on the {\it area\/} of a surface.  What is the correct statement of
the gluing law?  Of the symmetry properties?
        \endexer

 \subhead Beyond Quantum Hilbert Spaces
 \endsubhead

You should have noticed that we have not yet considered locality---a gluing
law---for the quantum Hilbert space.  This is not an idea which is usually
explicitly discussed in quantum field theory, though perhaps it is implicit
there.  For example, consider a discrete system formulated on a lattice in
space, i.e., on a $(d-1)$~dimensional lattice~$Y$.  Then the Hilbert space of
the theory is the tensor product over the lattice sites of a finite
dimensional Hilbert space~$H_y$ at each site:
  $$ E(Y) = \bigotimes\limits_{y\in Y}H_y.  $$
Then if the lattice~$Y$ is split into two pieces~$Y_1$ and~$Y_2$ the Hilbert
spaces obviously obey the equation
  $$ E(Y) \cong E(Y_1)\otimes E(Y_2). \tag{3.12} $$
This is what we mean by saying that the quantum Hilbert space is local,
though in general the gluing law is more complicated.  Also keep in mind what
we did for the classical theory.  The classical counterpart to the quantum
Hilbert space is the action~\thetag{2.2}, and the gluing law~\thetag{3.12} is
the quantum version of the (trivial case of the) gluing law~\thetag{2.16}.

The {\it Verlinde formulas\/}~\cite{V} are a nontrivial example of the
locality of the quantum Hilbert space.  It was originally formulated for the
spaces of ``conformal blocks'' in 2~dimensional conformal field theory.  Then
Witten~\cite{W} identified these spaces with the quantum Hilbert spaces of
3~dimensional Chern-Simons theory, and in this context the Verlinde formula
is a gluing law for the quantum Hilbert spaces.  It takes the following form.
Imagine that a closed surface~$Y$ is split into~$Y_1,Y_2$ along a circle~$S$.
There is a finite number of ``labels''~$\lambda  \in \Lambda $ and for each
label a Hilbert space~$E(Y_i)(\lambda )$.  The Verlinde formula roughly has
the form (we ignore some subtleties):
  $$ E(Y)\cong \bigoplus_{\lambda \in \Lambda }E(Y_1)(\lambda )\otimes
     E(Y_2)(\lambda ). \tag{3.13} $$

        \exertag{3.14} {Exercise}
 Reinterpret this formula along the following lines.  Introduce the
category~$\E$ whose objects are collections of Hilbert spaces indexed
by~$\Lambda $.  Also, introduce the ``inner product'' $\E\times \E\to\E$
defined by
  $$ \bigl(\{E_1(\lambda )\},\{E_2(\lambda )\} \bigr)_\E = \bigoplus_{\lambda
     }E_1(\lambda )\otimes E_2(\lambda ). \tag{3.15} $$
Show that this makes~$\E$ a 2-Hilbert space (cf.~\theexertag{2.17}
{Exercise}).  Now rewrite~\thetag{3.13} in terms of~$\E$.  Your formula
should look like~\thetag{3.8}.
        \endexer

So how should we prove a gluing law for the quantum Hilbert space?  The
easiest way would be to repeat the computation~\thetag{3.4}.  But that
requires that we write the quantum Hilbert space as an integral, as
in~\thetag{3.1}.  This is what we do!  It is really the crucial step in these
lectures.  Recall that our extended classical action~\thetag{2.2} takes
values in hermitian lines.  So we write exactly the same equation
as~\thetag{3.1}, replacing the classical action~$\eac XP$ by the classical
action~$\hl YQ$ for fields on a closed $(d-1)$-manifold:
  $$ E(Y) = \int_{\fld Y}\hl YQ\;\meas YQ. \tag{3.16} $$
What does this mean?  It is a finite sum
  $$ \mu _1\cdot L_1 +\dots + \mu _N\cdot L_N  $$
where the $\mu _i$~are positive numbers and the $L_i$~are hermitian lines.
Now we interpret $\mu \cdot L$~as the hermitian line with the same underlying
complex vector space as~$L$ but with the inner product multiplied by~$\mu $.
We interpret the sum as the orthogonal direct sum.  In this way
\thetag{3.16}~defines a Hilbert space.  In fact, it is exactly the same
Hilbert space as~\thetag{3.7}.

        \exertag{3.17} {Exercise}
 Verify this last assertion.  Recall that we use the $L^2$~inner product
on~\thetag{3.7}.
        \endexer

Now we repeat the computation~\thetag{3.4} for a $(d-1)$-manifold~$Y$ split
along a closed $(d-2)$-manifold~$S$.  In the course of that we will naturally
introduce
  $$ \E(S) = \text{space of sections of $\Cal{L}_S \to \fld{S}$}
      $$
which we could also write as an integral
  $$ \E(S) = \int_{\fld S}\Cal{L} _S(R)\;\meas SR. \tag{3.18} $$
We interpret this in the trivial case where $\Cal{L}_S(R)=\Cal{L}$, the
category of finite dimensional Hilbert spaces.  Then an element of~$\E(S)$ is
a choice of a Hilbert space~$W_R$ for each field~$R\in \fld S$.  Put
differently, an element of~$\E(S)$ is simply a hermitian vector bundle over
the space of fields~$\fld S$.  So
  $$ \E(S)=\Vect(\fld S) \tag{3.19} $$
is the collection of such hermitian vector bundles.

You can think of~$\E=\E(S)$ as being a ``Hilbert space over~$\LL$'', analogous
to the usual concept of a Hilbert space over~$\CC$.  We call such an object a
``2-Hilbert space''.  (2-vector spaces---the same object without the inner
product---were introduced by Kapranov/Voevodsky~\cite{KV} and
Lawrence~\cite{L}.) In \theexertag{2.11} {Exercise} we already made the
analogy between~$\CC$ and~$\LL$.  Now go back to \theexertag{2.14} {Exercise}
and \theexertag{2.17} {Exercise} to see what a 2-Hilbert space is.  There are
operations
  $$ \aligned
      +\:\E\times \E&\longrightarrow \E\\
      \cdot \:\LL\times \E &\longrightarrow \E\\
      (\cdot ,\cdot )\:\E\times \E&\longrightarrow \LL\endaligned \tag{3.20}
     $$
analogous to addition, scalar multiplication, and inner product.

        \exertag{3.21} {Exercise}
 Determine the operations~\thetag{3.20} for the 2-Hilbert space defined
in~\thetag{3.19}.  What is a basis for this 2-Hilbert space?
        \endexer

        \exertag{3.22} {Exercise}
 State explicitly the symmetry and gluing properties for the quantum
integrals~\thetag{3.1}, \thetag{3.16}, and~\thetag{3.18}.  Note in particular
the gluing law for the quantum Hilbert spaces:
  $$ E(Y) = \bigl(E(Y_1),E(Y_2) \bigr)_{\E(S)}.  $$
How does this fit with~\thetag{3.12}, \thetag{3.13}, and~\thetag{3.15}?
        \endexer

We summarize our extended notions of classical action and quantum path
integral in the following diagram.  Here all of the manifolds are presumed
closed for simplicity.

\def\mentry#1#2{$\vcenter{\vbox{\kern6pt\hsize%
 135pt\centerline{#1}\vskip12pt\centerline{#2}\kern6pt}}$}
 \bigskip
 \midinsert
 {\offinterlineskip \tabskip = 0pt
  \halign{
	\vrule\enspace\hfil#\hfil\enspace\vrule\hskip2pt
        \vrule&\enspace\hfil#\hfil\enspace
	&\vrule\enspace\hfil#\hfil\enspace\vrule\cr
 \noalign{\hrule}
 \titlestrut{\bf dimension}&{\bf classical action}&{\bf path integral} \cr
 \noalign{\hrule}
 \vphantom{\vrule height 2pt}&&\cr \noalign{\hrule}
 \entry $d\hphantom{-1.}\quad (X)$:\mentry{$P\longmapsto \eac
XP$}{complex number of unit norm}:\mentry{$Z_X$}{complex number}
 \noalign{\hrule}
 \entry $d-1\quad (Y)$:\mentry{$Q\longmapsto \hl YQ$}{hermitian
line}:\mentry{$E(Y)$}{Hilbert space}
 \noalign{\hrule}
 \entry $d-2\quad (S)$:\mentry{$R\longmapsto \LL_S(R)$}{one dimensional
2-Hilbert space}:\mentry{$\E(S)$}{2-Hilbert space}
 \vphantom{\vrule height2pt}&&\cr \noalign{\hrule}}
 }
 \nobreak
 \bigskip
 \centerline{Extended notions of classical action and path integral}
 \medskip
 \endinsert

 \subhead Quantum Finite Gauge Theory
 \endsubhead

We now compute the path integral in finite gauge theory, specializing
to~$d=3$.  In this section we treat the usual path integral and quantum
Hilbert space.  The new point is to account for the symmetry on the space of
fields when carrying out the quantization.  In Lecture~4 we will compute the
quantum object~$\E(\cir)$ associated to the circle and show how it leads to a
quantum group.

You should now review the last part of Lecture~2 to be sure you understand
the fields and classical action in this theory.  (Recall that we consider the
``untwisted'' case where the classical action is trivial.)  We resume that
discussion taking over the notation used there.  The ingredient we
are missing is a measure on the space of fields, which is simply defined:
  $$ \mu _M(P) = \frac{1}{\#\Aut P}, \qquad M\in \fld M.  \tag{3.23} $$
This is the correct ``counting measure'' and is always how we count objects
in mathematics.  Symmetries identify equivalent objects which we only want to
count once.  To be consistent, then, if an object has automorphisms it must
be counted as in~\thetag{3.23}.

        \exertag{3.24} {Exercise}
 Verify that this measure is invariant under symmetry.
        \endexer

The partition function is defined by~\thetag{3.1}, except now that we
integrate only over the space of {\it equivalence classes\/} of fields~$\fldb
X$.  This makes sense since both the classical action and the measure are
invariant under symmetries.  Writing~`$\Pbar$' for a typical equivalence
class we have
  $$ \split
      Z_X &= \int_{\fldb X}\eac X{\Pbar}\;\meas X{\Pbar}\\
      &=\sum\limits_{\Pbar}\frac{1}{\#\Aut P}.\endsplit \tag{3.25} $$

Suppose for simplicity that $X$~is connected.  Introduce a basepoint~$x\in X$
and consider the space of bundles with a basepoint:
  $$ \fld X' = \left\{ \langle\matrix P\\ \downarrow\\ X\endmatrix,p\rangle :
     P \text{ is a principal $G$~bundle, $p\in P_x$ a chosen basepoint
     .}\right\} \tag{3.26} $$
Here maps of bundles are required to preserve the basepoint.  We need a few
facts outlined in the next exercise.

        \exertag{3.27} {Exercise}
 Show that objects in~$\fld X'$ are {\it rigid\/}, that is, have no nontrivial
automorphisms.  In other words, once we specify an isomorphism of bundles on
the basepoint, we know it everywhere.  Show that the holonomy sets up a
1:1~correspondence between equivalence classes of elements in~$\fld X'$ and
homomorphisms $\pi _1(X,x)\to G$.  There is a $G$~action on~$\fld X'$: A
group element~$g\in G$ simply moves the pair~$\langle P,p \rangle$
to~$\langle P,p\cdot g \rangle$.  Show that this action passes to the
quotient~$\fld X'$.  Show that this quotient action corresponds to the action
of $G$ by conjugation on the space of homomorphisms $\Hom(\pi _1(X,x),G)$.
        \endexer

\flushpar
 Now the counting measure on~$\fld X'$ weights each bundle with basepoint
with weight~1, since there is only the identity automorphism.  Taking into
account the $G$~action we find
  $$ Z_X = \frac{\# \Hom(\pi _1(X),G)}{\#G}. \tag{3.28} $$
This is obviously a topological invariant of~$X$.

We now compute the quantum Hilbert space, defined by~\thetag{3.7}
or~\thetag{3.16}.  Let $Y$~be a closed surface.  To account for the gauge
symmetry we again integrate over the space of equivalence classes of fields.
What this means is that we replace~\thetag{3.7} by the space of {\it
invariant\/} sections, that is, sections invariant under gauge
transformations.  Since the action is trivial~\thetag{2.30}, this is merely
the space of invariant functions, i.e., the space of functions on the
quotient~$\fldb Y\cong \Hom\bigl(\pi _1(Y),G \bigr)/G$.  For this last
equality we assume that $Y$~is connected.  The $L^2$~metric is defined using
the measure~\thetag{3.23}.  Our argument with basepoints above identifies
this with
  $$ E(Y)\cong \frac{1}{\# G}\cdot L^2\bigl(\Hom(\pi _1(Y,y),G),G \bigr)^G.
     \tag{3.29} $$
Here $y$~is any basepoint in the connected space~$Y$, the symbol `$(\cdot )^G$'
means the invariants under the $G$~action by conjugation, the $L^2$~metric
weights each homomorphism with unit weight, and the prefactor~$1/\#G$
multiplies this $L^2$~metric.

        \exertag{3.30} {Exercise}
 Verify~\thetag{3.29}.  What is the answer for $Y=S^2$?  What about
$Y=\cir\times \cir$?  For this example note that $\Hom\bigl(\pi _1(Y,y),G
\bigr)$ is the set of commuting pairs of elements in~$G$.  To study the
action of conjugation on this set it helps to consider the picture in
\theexertag{2.18} {Exercise}.
        \endexer

        \exertag{3.31} {Exercise}
 By general arguments symmetries of~$Y$ are implemented as linear isometries
of~$E(Y)$.  Do some explicit computations for~$Y=\cir\times \cir$.  For
example, compute the effect of the diffeomorphism defined by the matrix
$T=\left(\smallmatrix 1&1\\0&1  \endsmallmatrix\right)\in SL(2;\ZZ)$.
        \endexer

        \exertag{3.32} {Exercise}
 Compute the path integral on a 3-manifold with boundary.  Verify the gluing
law for this relative invariant.
        \endexer

        \exertag{3.33} {Exercise}
 The $d=2$ case of this theory is also interesting.  It is the finite group
version of the zero area limit (or topological limit) of 2~dimensional QCD as
considered in the lectures of David Gross.  Now the partition function is
defined for surfaces and there is a basic quantum Hilbert space~$E$ attached
to the circle.  Show that $E$~can be identified as the space of {\it
central\/} functions on~$G$, that is, functions on~$G$ invariant under
conjugation.  What is the measure?  Show that the path integral over a ``pair
of pants'' leads to a multiplication on~$E$.  What is the multiplication?
Can you diagonalize it?  Use the gluing law to compute the partition function
on any closed surface.  Compare with~\thetag{3.28} to obtain a formula which
counts homomorphisms from a surface group into a finite group.
        \endexer

        \exertag{3.34} {Exercise}
 Try to do some computations for a twisted theory (in $d=3$) defined by a
nonzero~$\lambda \in H^3(BG;\RZ)$.  For example, take $G$~to be a cyclic
group of order~$n$.  Then $H^3(BG;\RZ)$ is also cyclic of order~$n$, and we
can take $\lambda $~to be the generator.  Compute the result for $X$~the
projective 3-space, or more generally a lens space.
        \endexer

\head
\S{4}: The Quantum Group
\endhead
\comment
lasteqno 4@ 40
\endcomment

Finally in this lecture we produce the quantum group in finite group
Chern-Simons theory.  For simplicity we only consider the untwisted theory,
and not the twisted theory which we indicated in the exercises of Lecture~2.
However, we should remark that the precise computations for the twisted
theory illustrate some more subtle gluing laws (along codimension~2
submanifolds).  The interested reader should consult~\cite{F1,\S\S8--9} for
details.

Our arguments in this section are somewhat sketchy.  A rigorous treatment
would at the very least demand that we make precise all of the axioms for
2-Hilbert spaces, and this is already a complicated matter.  The reader may
refer to~\cite{F1,\S5,\S7} for another treatment.

 \subhead The 2-Hilbert Space
 \endsubhead

We now compute the 2-Hilbert space~$\E=\E(\cir)$ attached to the standard
oriented circle.  (We take $\cir$~to be the unit circle in the complex
numbers with the counterclockwise orientation.  We also take $1\in \cir$ as a
basepoint.)  Recall from~\thetag{3.18} that this is defined as a
(generalized) integral over the space of fields on the circle.  So our first
job is to understand what that space $\fld{\cir}$ looks like.  As
in~\thetag{3.26} it is best to start by rigidifying the fields by introducing
a basepoint.  So we consider $\fld{\cir}'$, the category of principal
$G$~bundles over a circle with a chosen basepoint covering the basepoint
in~$\cir$.  (For brevity we call them ``pointed principal bundles.'')  Then
by starting at the basepoint~$p$ and traversing the circle in the properly
oriented direction, we arrive back at a point~$p\cdot x$, where $x\in G$~is
the {\it holonomy\/}.  (See Figure~5.)  So the holonomy around~$\cir$ defines
a map
  $$ \hol\:\fld{\cir}'\longrightarrow G,  $$
and it is easy to verify that it is an isomorphism on equivalence classes:
  $$ \overline{\fld{\cir}'}\cong G.  $$
Now a change of basepoint $p\to p\cdot g$ changes the holonomy from~$x$
to~$g\inv xg$.  Thus the conjugacy class of the holonomy is independent of
the basepoint, and
  $$ \fldb{\cir} \cong \text{conjugacy classes in~$G$}, \tag{4.1} $$
which agrees with~\thetag{2.28}.

We know by now that it is not enough to work with equivalence classes.  So as
at the end of the last lecture we use basepoints to construct a good model
of~$\fld{\cir}$.  Namely, introduce the category~$\C$ whose objects are
elements~$x\in G$ and with a morphism $x @>g>> gxg\inv $ for every pair of
elements~$\langle x,g \rangle$.  (See \theexertag{2.18} {Exercise}.)  Think
of the object labeled by~$x$ as a pointed principal bundle~$P_x$ with
holonomy~$x$, and the morphism labeled by $x@>g>>gxg\inv $ as the map between
the pointed principal bundles $P\mstrut _x\to P_{gxg\inv }$ which takes the
basepoint in~$P_x$ to $g$~times the basepoint in~$P_{gxg\inv }$.  These
morphisms are isomorphisms of (unpointed) principal bundles, and represent
all of the possible morphisms.  For example, the set of arrows $x@>g>>x$ is
isomorphic to the {\it centralizer\/}~$C_x$ of~$x$ in~$G$.  This is precisely
the set of automorphisms of~$P_x$ as a principal bundle.  Such arguments show
that $\C$~is a good model of~$\fld{\cir}$.

Of course, $\C$~is simply a picture of the action of~$G$ on itself by
conjugation.

We can also identify an arrow $\arr xg gxg\inv $ as a field on the
cylinder~$C$.  We picture that in Figure~6, where we assume basepoints at
the top of each circle and the group element indicates the parallel transport
relative to the basepoints.  (The cylinder is depicted as an annulus.)

\midinsert
\bigskip
\centerline{
\epsffile{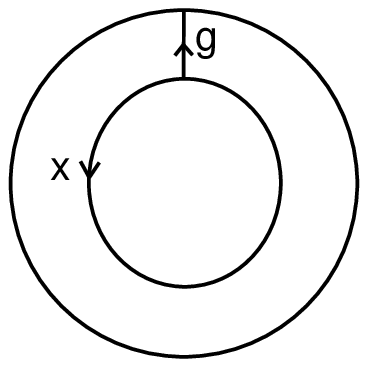}}
\nobreak
\medskip
\centerline{Figure~6: The bundle on~$C$ corresponding to $\arr xg$}
\bigskip
\endinsert

Let $\G$~denote the set of arrows in~$\C$; it is a {\it groupoid\/}.  (A
groupoid is a semigroup in which every arrow is invertible.)  For simplicity
we denote elements of~$\G$ by~$\arr xg$.  Then the composition law in~$\G$ is
  $$ (\arr {x'}{g'})\;\circ\;(\arr xg ) = \cases \arr x{g'g} \;,&x'=gxg\inv
     ;\\\text{undefined},&x'\not= gxg\inv .\endcases \tag{4.2} $$
We can understand this by gluing two bundles on the cylinder.  The inner
bundle should correspond to~$\arr xg$ in keeping with our convention of
representing parallel transport as acting on the left.  (It is, after all, a
bundle morphism.)

Recall from~\thetag{2.31} that in this theory the classical action is
trivial.  (The twisted theory has a nontrivial action.)  Now if there were no
symmetries and the space of fields were simply~$\fldb{\cir}$, which
by~\thetag{4.1} is the set of conjugacy classes in~$G$, then by~\thetag{3.19}
we would identify~$\E$ with the category of hermitian vector bundles over the
set of conjugacy classes in~$G$.  In other words, an element of~$\E$ would be
a collection of hermitian vector spaces indexed by the conjugacy classes
of~$G$.  However, this is not correct because of the nontrivial symmetries in
the fields.

When we quantize in the top dimension we deal with the symmetries by
performing the path integral over the space of equivalence classes of
fields~\thetag{3.25}.  In the next dimension down, when we quantize the
theory on a surface, we accommodate the symmetries by considering {\it
invariant\/} sections~\thetag{3.29}.  Now we must implement the analogous
idea in codimension two: We take ``invariant sections'' of the trivial bundle
whose fiber at each point is~$\catl$, the category of finite dimensional
Hilbert spaces.  But rather than consider the large category~$\fld{\cir}$, we
use our model~$\C$ for the space of fields on the circle.  Then a section is
simply a hermitian vector bundle $W\to G$.  In other words, it is a
collection of hermitian vector spaces~$W_x$ indexed by~$x\in G$.  What do we
mean by an ``invariant'' section?  Well, we mean that it should be invariant
under the morphisms in~$\C$, that is under the action of~$G$ on itself by
conjugation.  However, we do not simply mean that $W_x=W_{gxg\inv }$.
Rather, we mean that we are given an explicit isomorphism
  $$ A_g\:W\mstrut _x\longrightarrow W_{gxg\inv } \tag{4.3} $$
for each $x,g\in G$.  This is in line with the idea that sets have an
additional layer of structure: Vector spaces can be isomorphic without being
equal.  Of course, we presume that the isomorphisms~\thetag{4.3} compose in
accordance with the composition law~\thetag{4.2} for the arrows.  So an
element of~$\E$ is simply a vector bundle~$W\to G$ together with a lift of
the conjugation action of~$G$ on itself.  The collection of these {\it
equivariant vector bundles\/} is usually denoted~$\Vect_G(G)$.  Observe that
the dimension of the fiber~$W_x$, though constant in a conjugacy class, can
vary over the group~$G$.  Also, the bundles we consider have a {\it
hermitian\/} structure.

There is one other ingredient in the definition~\thetag{3.18} of~$\E$, namely
the measure on the space of fields.  The correct measure to use on unpointed
bundles is~\thetag{3.23}.  However, our pointed bundles are rigid; they don't
have any nontrivial automorphisms.  But since we have to ``divide'' by the
action of~$G$ (by taking invariant sections) we should use the weight
$1/\#G$.  (Precisely the same factor occurs when we quantize a closed
surface~\thetag{3.29}.)  So finally,
  $$ \E \cong \frac{1}{\#G}\cdot \Vect_G(G).  $$
The standard inner product in~$\catl$ is $(V_1,V_2) = V_1 \otimes
\overline{V_2}$, and that in~$\Vect(G)$ is obtained by summing over the
fibers.  To account for the $G$~action we need to take the invariants.
Putting this together we see that the inner product in~$\E$ is
  $$ (W_1,W_2)_{\E} = \frac{1}{\#G}\cdot \( \bigoplus_x (W_1)_x\otimes
     \overline{(W_2)_x}\) ^{\tsize G} . \tag{4.4} $$
Note that the right hand side is a finite dimensional Hilbert space, as it
should be.

        \exertag{4.5} {Exercise}
 Make explicit the 2-Hilbert space structure of~$\E$.  What is the addition?
Scalar multiplication?
        \endexer

        \exertag{4.6} {Exercise}
 Rewrite~\thetag{4.4} as a sum over conjugacy classes.
        \endexer

        \exertag{4.7} {Exercise}
 Compute~$\E$ for an abelian group~$G$.

        \endexer

Next, we observe that $\E$~can be identified as the category of
representations of the groupoid~$\G$.  (In \theexertag{2.17} {Exercise} we
indicate how the category of representations of a finite {\it group\/} is a
2-Hilbert space.  Here we replace a finite group by a finite groupoid.)
Namely, associate to an element $W\in \Vect_G(G)$ the finite dimensional
Hilbert space
  $$ W = \bigoplus_x W_x.  $$
(The overloaded notation should not cause confusion.)  Then the arrow $(\arr
xg)$ acts on~$W$ by~$A_g\res{W_x}$; the action is trivial on $W_{x'}$
for~$x'\not= x$.  We can recover the fibers of the vector bundle by $W_x =
(\arr xe)(W)$, where $e$~is the identity element of~$G$.

We reformulate this by introducing an algebra~$H$ which we might call the
``groupoid algebra,'' by analogy with group algebras:
  $$ H = \bigoplus_{x,g}\,\CC\langle x,g \rangle.  $$
(We use the symbol `$\langle x,g \rangle$' for an element of~$H$ to
distinguish it from the corresponding element $(\arr xg)$ in~$\G$.)  So an
element of~$H$ is a formal linear combination of the symbols~$\langle x,g
\rangle$ with complex coefficients.  The multiplication in~$H$ is
  $$ \langle x_2,g_2 \rangle\cdot \langle x_1,g_1 \rangle = \cases \langle
     x_1,g_2g_1 \rangle ,&x\mstrut _2=g\mstrut _1x\mstrut _1g_1\inv
     ;\\0,&\text{otherwise} .\endcases \tag{4.8} $$
The unit element is
  $$ 1=\sum\limits_{x}\langle x,e \rangle.  \tag{4.9} $$
Then $\E$~is the 2-Hilbert space of unitary representations of the
algebra~$H$, except that the usual inner product is divided by~$\#G$.

We pass freely among these various descriptions of~$\E$.

What is a basis for~$\E$, thought of as a 2-Hilbert space?  The geometric
picture of equivariant vector bundles~$W\to G$ is perhaps easiest to
consider.  Note that the fiber~$W_x$ is a representation of the
centralizer~$C_x$, and that the representations~$W_{x'}$ and~$W_x$ are
isomorphic if $x$~and $x'$~are conjugate.  So an indecomposable element
of~$\Vect_G(G)$ is supported on a single conjugacy class and for each~$x$ in
that conjugacy class is an irreducible representation of the
centralizer~$C_x$.  Let $\{W_\lambda \}$~be a set of inequivalent irreducible
elements.  Then it is a basis for~$\E$ in the sense that any other element is
isomorphic to $\oplus_{\lambda }\, (V^\lambda \otimes W_\lambda) $ for
some finite dimensional Hilbert spaces~$V^\lambda $ (thought of as trivial
vector bundles over~$G$).

        \exertag{4.10} {Exercise}
 Write a basis of~$\E$ for $G$~an abelian group.  For $G=S_3$, the symmetric
group on 3~letters.
        \endexer

 \subhead Locality and Gluing
 \endsubhead

So far we have just used the extended notion of quantization and the
structure of the space of fields on the circle.  Now we want to apply the
basic principles of locality and gluing to derive further structure on~$\E$.
At each stage this can be realized by introducing further structure on~$H$.
By then end we will have introduced enough structure to make~$H$ into a
quantum group.  We begin by considering the path integral over surfaces with
boundary.

If $Y$~is any surface with boundary, then from our general picture in
Lecture~3 we know that the path integral~$E(Y)$ is an element of~$\E(\bY)$.
Now $\bY$~is a union of circles, and we would like to identify each of these
circles with the standard circle so that $\E(\bY)$~can be identified as a
tensor product of copies of~$\E=\E(\cir)$.  In general we cannot do this
without introducing additional data (such as a parametrization of the
boundary circles), but we will only need to consider surfaces~$Y\subset \CC$,
that is, disks with subdisks removed.  Then the boundary is made up of
circles in~$\CC$, and there is a unique composition of a translation and a
dilation which brings each circle to the standard circle.  In this way we can
identify $E(Y)$ as living in a tensor product of copies of~$\E$.  We can
designate some of the boundary components as `incoming' and some as
`outgoing' and use the hermitian metric in~$\E$ to view
  $$ E(Y) \: \bigotimes_{\text{incoming circles}} \!\!\!\!\!\!\!\!\!\!\!\E
     \;\;\;\longrightarrow \bigotimes_{\text{outgoing
     circles}}\!\!\!\!\!\!\!\!\!\!\!\E .   $$
This is a ``linear'' map of 2-Hilbert spaces, which at the simplest level is
a functor between the categories indicated.

Introduce a basepoint on each boundary component of~$Y$.  Let $\fld Y'$
denote the category of principal $G$~bundles over~$Y$ with chosen basepoints
over the basepoints in~$\bY$.  Assume that $Y$~has no closed components.
Then arguments similar to those at the end of Lecture~3 show that we can
identify
  $$ E(Y)\cong L^2(\overline{\fld Y')}. \tag{4.11} $$
Here we think of~$\E$ as the 2-Hilbert space of representations of the
groupoid~$\G$, and the $\G$~action of a given boundary component is computed
by gluing a cylinder to that boundary component.

        \exertag{4.12} {Exercise}
 Verify~\thetag{4.11} in detail.
        \endexer

        \exertag{4.13} {Exercise}
 Suppose that $f\:Y'\to Y$ is a diffeomorphism, where $Y,Y'\subset \CC$ are
surfaces as in the previous exercise.  Explain how $f$~induces an
automorphism of functors~$E(Y')\to E(Y)$.
        \endexer

        \exertag{4.14} {Exercise}
 Let $C$~denote the cylinder.  Compute~$E(C)$.
        \endexer

We will compute the path integrals over some surfaces shortly, but we begin
by considering symmetries of the circle.

        \exertag{4.15} {Exercise}
 Consider the group of orientation-preserving diffeomorphisms~$\Diff^+(\cir)$
of the circle.  Inside this group is the subgroup of rigid rotations.
Construct a retraction of~$\Diff^+(\cir)$ onto this subgroup.  (Hint:  Given
a diffeomorphism compose with a rotation so that the composition preserves
the basepoint.  Now cut open the circle at the basepoint; then an
orientation-preserving diffeomorphism is a monotonically increasing function
$f\:\zo\to\zo$, and if it preserves the basepoint then $f(0)=0$ and
$f(1)=1$.)
        \endexer

Recall from \theexertag{3.9} {Exercise} that diffeomorphisms of surfaces
induce isometries of the corresponding quantum Hilbert spaces.  In
\theexertag{3.22} {Exercise} you were asked to state the corresponding
symmetry, or functoriality, for the 2-Hilbert space associated to a
1-manifold (or a $(d-2)$-manifold in a $d$~dimensional theory).  At first
glance we might expect a diffeomorphism $f\:S'\to S$ to induce an isometry of
2-Hilbert spaces $f_*\:\E(S')\to \E(S)$, and indeed this is true.  However,
again there is an extra layer of structure: We expect an isotopy~$F=f_t$
(path of diffeomorphisms) from $f_0$ to~$f_1$ to induce a natural
transformation~$F_*$ from the functor~$(f_0)_*$ to the functor~$(f_1)_*$.  In
topological theories we expect that deformations will not change this natural
transformation.  Applied to $f_0=f_1=id_S$, we obtain an action of~$\pi
_1(\Diff^+(S))$ by automorphisms of the identity functor.  Such an {\it
automorphism of the identity\/} commutes with all arrows in the
category~$\E$.  More explicitly, an automorphism of the identity specifies
for each~$W\in \E $ a map
  $$ \theta _W\:W\longrightarrow W  $$
which commutes with all arrows in~$\E$.

        \exertag{4.16} {Exercise}
 Let $\Rep(H)$ denote the 2-Hilbert space of representations of a finite
group~$H$. (See \theexertag{2.17} {Exercise}.)  Determine all automorphisms
of the identity of~$\Rep(H)$.
        \endexer

We apply this to the circle.  From the previous exercise we know that
$\Diff^+(\cir)$ is homotopy equivalent to a circle, so the fundamental group
is infinite cyclic.  It is generated by a circle of rotations.  Write an
element of the circle as~$e^{2\pi ix}$; then the circle of rotations is
  $$ f_t(e^{2\pi ix}) = e^{2\pi i(t+x)},\qquad 0\le t\le 1.  $$
We view this circle of rotations as a diffeomorphism~$\tau $ of the
cylinder~$C$ which is the identity on both boundary circles.  We compute the
corresponding automorphism of the identity by computing the action on~$E(C)$.
(See \theexertag{4.13} {Exercise}.)  Now since the cylinder glued to itself
is the cylinder, it follows easily that $E(C)$~is the identity map~$\E\to\E$
(\theexertag{4.14} {Exercise}).  We first compute the induced action on
fields on the cylinder as
  $$ \tau ^*\langle x,g  \rangle = \langle x,gx  \rangle = \langle x,g
     \rangle\cdot \langle x,x  \rangle. \tag{4.17} $$
Here we use the multiplication in the algebra~$H$~\thetag{4.8}.  So the
action on a field with holonomy~$x$ may be described by gluing on the
field~$\langle x,x \rangle$.  Then \theexertag{4.12} {Exercise} implies that
the effect on the quantization is the operator $A_x\:W_x\to W_x$, i.e., the
automorphism of the identity is
  $$ \theta _W\res{W_x} = A_x\:W_x\longrightarrow W_x.   $$
In terms of the $H$~action we can describe it as the action of the special
element
  $$ v=\sum\limits_{x}\langle x,x  \rangle. \tag{4.18} $$
This special element is the inverse of the {\it ribbon element\/} defined by
Reshetikhin and Turaev~\cite{RT2}.

        \exertag{4.19} {Exercise}
 Show that \thetag{4.17}~induces a diffeomorphism of the torus.  In fact, it
is the diffeomorphism considered in \theexertag{3.31} {Exercise}.  Compare
the automorphism of the identity computed here with the action of that
diffeomorphism on the vector space associated to the torus.  Is there a
gluing law which compares them?
        \endexer

Purely abstract considerations (Schur's lemma) show that an automorphism of
the identity of~$\E$ acts as multiplication by a scalar on an irreducible
element~$W_\lambda $.  We can also see that since $x$~is a central element of
the centralizer group~$C_x$, so acts as a scalar in any irreducible
representation.  These scalars are the {\it conformal weights\/} of the
theory.

        \exertag{4.20} {Exercise}
 Fix an element~$x$ and a representation~$\lambda $ of the centralizer~$C_x$.
This determines a basis element of~$\E$.  Compute the conformal weight in
terms of the character of~$\lambda $.
        \endexer

        \exertag{4.21} {Exercise}
 Calculate the conformal weights explicitly for a cyclic group.  For the
symmetric group~$S_3$.
        \endexer

The circle also has an {\it orientation-reversing\/} diffeomorphism---a
reflection---which is unique up to isotopy.  We expect it to induce a map
$\E\to\overline{\E}$ whose square is the identity (since a reflection squares
to the identity).  Now $\overline{\E}$~means the conjugate 2-Hilbert space,
but as a map of categories we can ignore the conjugation---that only affects
the ``scalar'' multiplication.  Again we can compute the induced action
on~$\E$ from the cylinder.  Let $\rho $~denote a fixed reflection.  Then the
induced map on fields is
  $$ \rho ^*\langle x,g  \rangle = \langle x\inv ,g  \rangle.  $$
Note however that the incoming and outgoing circles are exchanged.  This
means that an element $W\in \E$ is taken to a new element~$W^*$ with
  $$ \aligned
      (W^*)\mstrut _x&=W_{x\inv }^*\\
      A^{W^*}_g &= (A^W_{g\inv })^*.\endaligned  $$
In terms of the algebra~$H$ the duality is implemented by the {\it
antipode\/}
  $$ S(\langle x,g  \rangle) = \langle gx\inv g\inv ,g\inv   \rangle.
     \tag{4.22} $$

        \exertag{4.23} {Exercise}
 In the diffeomorphism group of the circle, write the equation which relates
a rotation and a reflection.  In the quantization this is reflected by a
relation among conformal weights.  Namely, show that
  $$ \theta \mstrut _{W^*} = \theta _W^*.   $$
What does this say about conformal weights in terms of our basis?
        \endexer

Next, we consider the path integral over some simple surfaces~$Y\subset \CC$.
For example, consider the unit disk~$D$.  Since the disk is contractible, all
$G$~bundles over~$D$ are trivial, as any two basepoints are related by an
automorphism we find $\overline{\fld D'}$~consists of a unique element.  It
follows that $E(D)\cong L^2(\overline{\fld D'})$~is one dimensional and is
described by
  $$ E(D)_e = \cases \CC ,&x=e;\\0,&x\not= e,\endcases  \tag{4.24} $$
with $C_e=\Gamma $~acting trivially on~$E(D)_e$.  In terms of the algebra~$H$
this particular representation is described by a {\it counit\/} $\epsilon
\:H\to\CC$.  From~\thetag{4.24} we see that the counit is
  $$ \epsilon \bigl(\langle x,g \rangle\bigr) = \cases 1
     ,&x=e;\\0,&\text{otherwise} .\endcases \tag{4.25} $$

\midinsert
\bigskip
\centerline{
\epsffile{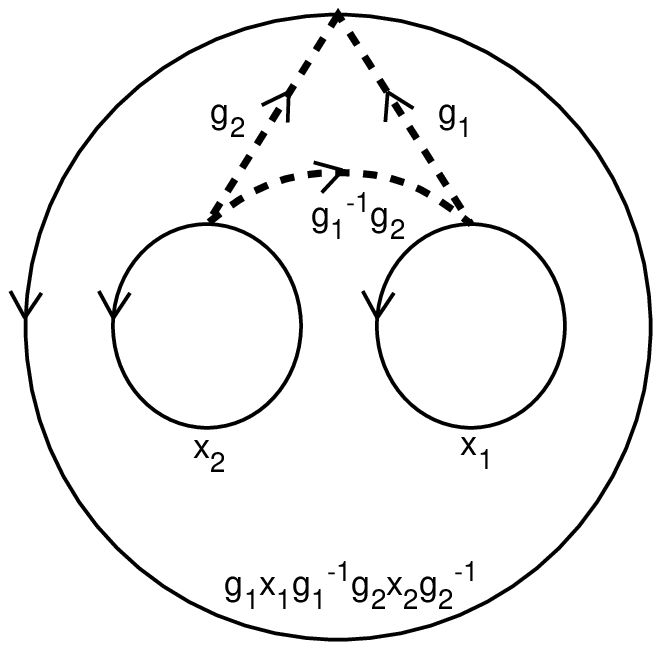}}
\nobreak
\medskip
\centerline{Figure~7: The bundle over~$P$ corresponding to $\langle x_1,g_1
\rangle\times \langle x_2,g_2  \rangle\in \G\times \G$}
\bigskip
\endinsert

More interesting is the path integral over the ``pair of pants,'' which we
realize embedded in~$\CC$ as the surface~$P$ depicted in Figure~7.  We view
the inside circles as incoming and the outside circle as outgoing (note the
orientations), so that $E(P)\in \E\otimes \E\otimes \E$ induces a map
  $$ \mytimes:\ \E\otimes \E\longrightarrow \E.  $$
Here we denote the result of this multiplication map on~$W_1,W_2$
as~$W_1\mytimes W_2$.  There is an induced map on morphisms in~$W_i$ as well:
$\mytimes$~is a functor!  (In category language it gives a ``monoidal''
structure to~$\E$.)

        \exertag{4.26} {Exercise}
 Use the gluing law to describe an associativity isometry
  $$ \varphi _{W_1,W_2,W_3}\:(W_1\mytimes W_2)\mytimes W_3\longrightarrow
     W_1\mytimes (W_2\mytimes W_3).  $$
(See Figure~8.)  In fact, you will find it to be trivial.  However, in other
theories (e.g. the twisted version of this theory) it is nontrivial.
        \endexer

        \exertag{4.27} {Exercise}
 Use gluing to show that multiplication by \thetag{4.24}~is isomorphic to the
identity map.  How is this expressed precisely?
        \endexer

\midinsert
\bigskip
\centerline{
 \epsfxsize=350pt
\epsffile{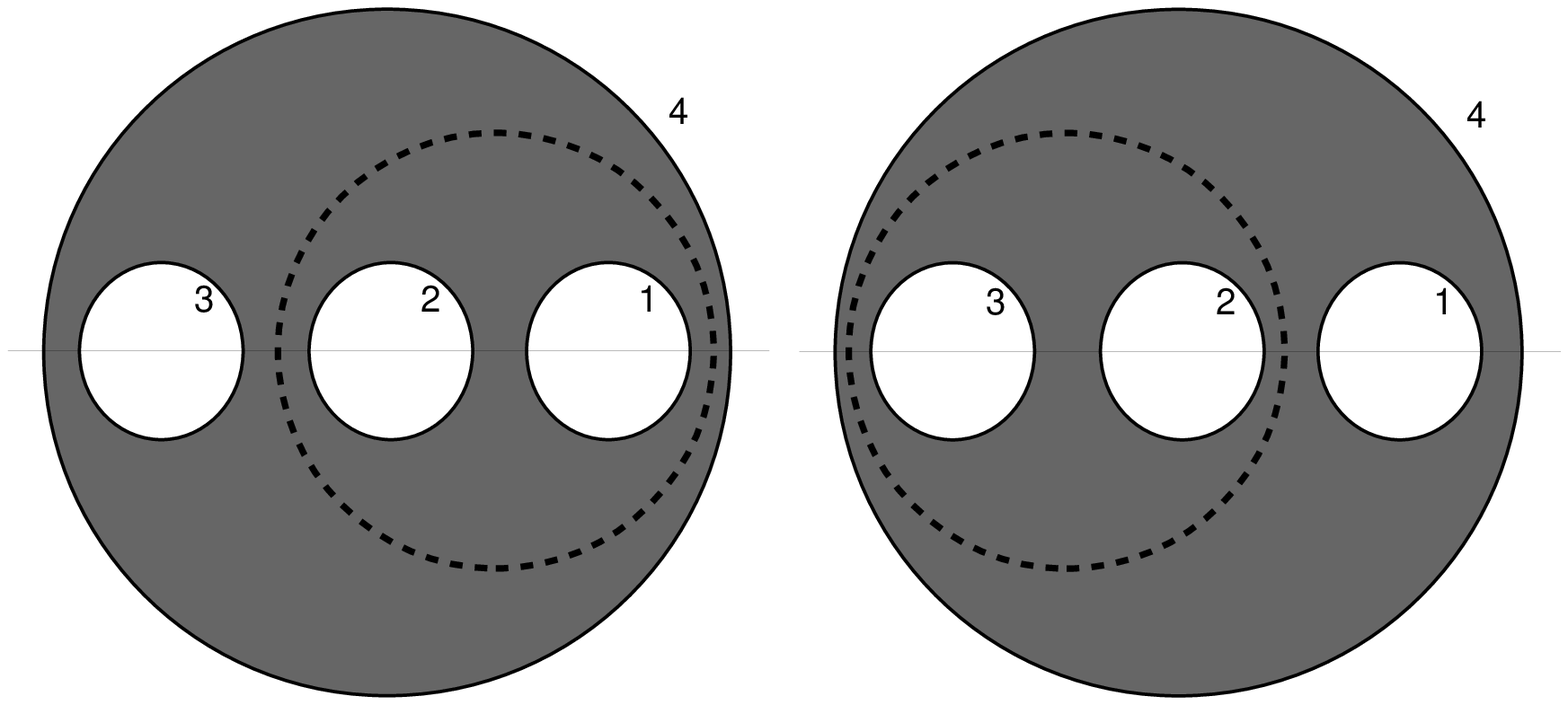}}
\nobreak
\medskip
\centerline{Figure~8: Associativity}
\bigskip
\endinsert

To compute~$E(P)$ we need to determine the space of fields on~$P$.  Following
our procedure, we introduce basepoints on each boundary component.

        \exertag{4.28} {Exercise}
 Show that there is a 1:1 correspondence between~$\fldbp P$ and $\G\times
\G$.  Hint:  See Figure~7.  In that figure the group elements indicate the
parallel transport along the given line relative to the chosen basepoints.
        \endexer

Now we must compute the three different $\G$~actions on~$E(P)\cong L^2(\fldbp
P)$ corresponding to the three boundary circles of~$P$.  For the inner two
boundary circles it is not difficult to verify that the action is by right
multiplication.  Namely, the action of~$\langle x,g  \rangle\in \G$ on a
function~$f(\cdot ,\cdot )$ on~$\fldbp P\cong \G\times \G$ is:
  $$ \aligned
     f(\cdot ,\cdot )&\longmapsto f(\cdot \langle x,g  \rangle\inv ,\cdot ),\\
     f(\cdot ,\cdot )&\longmapsto f(\cdot,\cdot\langle x,g  \rangle\inv
     ).\endaligned  \tag{4.29} $$
On the other hand, the action corresponding to the outer component is
  $$ \bigl(\langle x,g \rangle\cdot f \bigr)\bigl(\langle x_1,g_1
     \rangle,\langle x_2,g_2 \rangle \bigr) = \cases f\bigl(\langle x_1,gg_1
     \rangle,\langle x_2,gg_2 \rangle \bigr) ,&\text{if $x = g\mstrut
     _1x\mstrut _1g_1\inv g\mstrut _2x\mstrut _2g_2\inv$}
     ;\\0,&\text{otherwise} .\endcases \tag{4.30} $$

        \exertag{4.31} {Exercise}
 Verify \thetag{4.29} and \thetag{4.30}.  Recall that in each case the action
is induced by gluing on a cylinder to the appropriate boundary component.
        \endexer

To compute $W_1\mytimes W_2$, we use the inner product~\thetag{4.4} in~$\E$
to contract the first two factors (corresponding to the inner circles)
in~$E(P)\in \E\otimes \E\otimes \E$ with~$W_1\otimes W_2$.  Up to a factor in
the inner product, this says that $W_1\mytimes W_2$ is the $\G\times
\G$~invariants in~$W_1\otimes W_2\otimes E(P)$.  It is not difficult to
verify that
  $$ (W_1\mytimes W_2)_x = \bigoplus_{x_1x_2=x}(W_1)_{x_1}\otimes (W_2)_{x_2}
     \tag{4.32} $$
with $G$~action given by
  $$ A^{W_1\mytimes W_2}_g = A^{W_1}_g \otimes A^{W_2}_g. \tag{4.33} $$

        \exertag{4.34} {Exercise}
 Verify \thetag{4.32}--\thetag{4.33}.  Also, show that the
multiplication~$\mytimes$ can be described by the following operation on
{}~$\Vect_G(G)$.  Namely, let $\mu \:G\times G\to G$ denote multiplication in
the group~$G$.  Then given $W_1,W_2\in \Vect_G(G)$, we can form their
external tensor product ~$W_1\boxtimes W_2\to G\times G$.  The quantum
product is
  $$ W_1\mytimes W_2 = \mu _*(W_1\boxtimes W_2),  $$
where $\mu _*$ is the pushforward of vector bundles.  (This is defined for
finite covering maps.)
        \endexer

Viewing elements of~$\E$ as representations of the algebra~$H$, a tensor
product on representations is induced by a {\it coproduct\/} $\Delta \:H\to
H\otimes H$, since then the element $h\in H$ acts in the tensor
product~$W_1\otimes W_2$ by~$\Delta (h)$.  From~\thetag{4.32}
and~\thetag{4.33} we see that the appropriate coproduct is
  $$ \Delta \bigl(\langle x,g \rangle\bigr) = \sum\limits_{x_1x_2=x}\langle
     x_1,g \rangle\otimes \langle x_2,g \rangle.  \tag{4.35} $$

We content ourselves with one more computation, that of the {\it
$R$-matrix\/}.  Namely the pair of pants~$P$ has a {\it braiding\/}
diffeomorphism~$\beta $ which acts on the quantization~$E(P)$.  We depict the
braiding in Figure~9.  The diffeomorphism~$\beta $ rotates circle~1 under
circle~2, exchanging their positions.  It is the identity on circle~3.  Now
by \theexertag{4.13} {Exercise} $\beta $~induces an automorphism of the
multiplication~$\mytimes$.  Because two of the boundary circles are switched,
this means that for each~$W_1,W_2\in \E$ we have a map
  $$ R_{W_1,W_2}\:W_1\mytimes W_2\longrightarrow W_2\mytimes W_1
     $$
which commutes with all of the arrows in~$\E$.  Before proceeding to compute
it, we note one property which follows from the structure of~$\Diff^+(P)$.

\midinsert
\bigskip
\centerline{
  \epsfxsize= 350pt
\epsffile{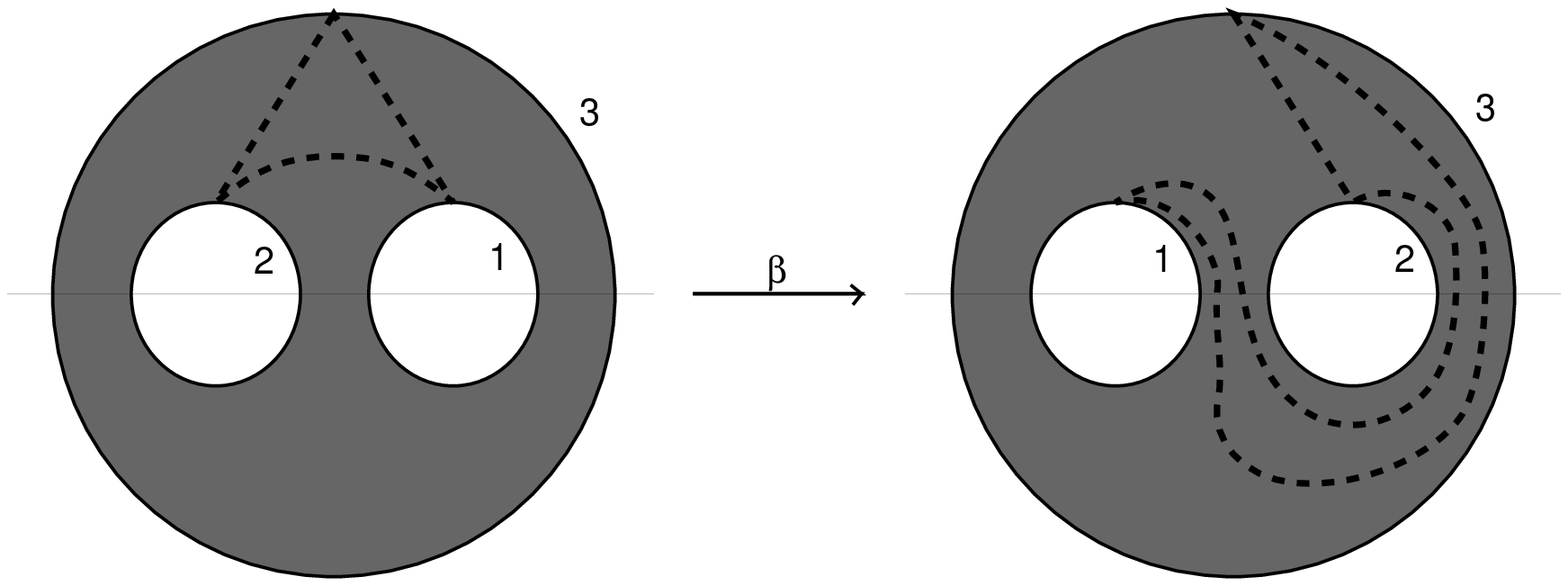}}
\nobreak
\medskip
\centerline{Figure~9: The braiding diffeomorphism~$\beta $}
\bigskip
\endinsert

        \exertag{4.36} {Exercise}
 Use pictures like Figure~9 to prove that
  $$ \tau _2 \tau _1 \beta = \beta \inv  \tau _3,  $$
where $\tau _i$~is a positive Dehn twist around the boundary labeled~$i$.
Show that this implies that for all~$W_1,W_2\in \E$ the following diagram
commutes:
  $$ \CD
      W_1\mytimes W_2 @>R_{W_1,W_2}>> W_2\mytimes W_1\\ @V{\theta _{W_1\mytimes
     W_2}}VV @VV{\theta _{W_2}\mytimes \theta _{W_1}}V\\
      W_1\mytimes W_2 @>R_{W_2,W_1}\inv >> W_2\mytimes W_1\endCD
     $$
        \endexer

Now to the computation of~$R_{W_1,W_2}$.  As in all computations the idea is
to first see the induced action on the fields.  This is easily computed to be
  $$ \langle x_1,g_1 \rangle\times \langle x_2,g_2 \rangle@>{\hbox to
     25pt{\hfil $\beta _*$\hfil}}>> \langle x_2\mstrut ,g\mstrut _1x\mstrut
     _1g_1\inv g\mstrut _2 \rangle\times \langle x_1,g_1 \rangle.
      $$
Thus on functions~$f\in L^2(\fldbp P)\cong E(P)$ the induced pushforward
action is:
  $$ \bigl(\beta _*f \bigr)\bigl(\langle x_1,g_1 \rangle,\langle x_2,g_2
     \rangle \bigr) = f\bigl(\langle x\mstrut _2,g\mstrut _1x\mstrut
     _1g_1\inv g\mstrut _2 \rangle,\langle x\mstrut _1,g\mstrut _1 \rangle
     \bigr).  $$
Now we must translate this into an action on~$W_1\mytimes W_2$ using the
construction preceding~\thetag{4.32}.  The result is
  $$ \aligned
       R_{W_1,W_2}\: (W_1)\mstrut _{x_1}\otimes (W_2)\mstrut
     _{x_2}&\longrightarrow (W_2)\mstrut _{x\mstrut _1x\mstrut _2x_1\inv
     }\otimes (W_1)\mstrut _{x_1} \\
       w_1\otimes w_2 &\longmapsto A_{x_1}^{W_2}(w_2)\otimes w_1\endaligned
     \tag{4.37} $$

        \exertag{4.38} {Exercise}
 Verify~\thetag{4.37}.
        \endexer

Finally, this can be implemented universally by a {\it quasitriangular
element\/}~$R\in H\otimes H$ which satisfies
  $$ R_{W_1,W_2} = \tau _{W_1,W_2}\circ (\rho _1\otimes \rho _2)(R),
      $$
where $\tau _{W_1,W_2}\:W_1\otimes W_2\to W_2\otimes W_1$ is the
transposition and $\rho _i$~are representations of~$H$. From~\thetag{4.37} we
deduce
  $$ R = \sum\limits_{x_1,x_2} \langle x_1,e \rangle\otimes \langle x_2,x_1
     \rangle.   \tag{4.39} $$

Thus we have arrived at our goal.  Namely from the quantization~$\E$ of the
circle, using symmetry and gluing laws, we have constructed an algebra~$H$
with unit~\thetag{4.9}, counit~\thetag{4.25}, antipode~\thetag{4.22},
comultiplication~\thetag{4.35}, a quasitriangular element~\thetag{4.39}, and
a ribbon element~\thetag{4.18}.  This is a {\it quasitriangular Hopf
algebra\/} with a ribbon element.  It is certainly the proper realization of
a {\it quantum group\/} in this example.  It was first written down in a
paper of Dijkgraaf, Pasquier, and Roche~\cite{DPR}, and it can be identified
with Drinfeld's {\it quantum double\/} of the Hopf algebra of functions
on~$G$.

        \exertag{4.40} {Exercise}
 Compute the Hopf algebra for $G$~a cyclic group.  For~$G=S_3$.
        \endexer

\enddocument